\documentclass{article}

\usepackage{hyperref}

\usepackage{lineno}
\usepackage[margin=1in]{geometry}
\usepackage{float}
\usepackage[]{algorithm}
\usepackage{algpseudocode}
\usepackage{graphicx}
\usepackage{caption}

\usepackage{subcaption}
\usepackage{longtable}
\usepackage{xcolor}
\usepackage{wrapfig}
\usepackage{tabularx}
\usepackage{authblk}
\usepackage[title, titletoc]{appendix}

\usepackage{amssymb,amsmath,amsthm}

\usepackage{booktabs}
 \allowdisplaybreaks
\usepackage{stmaryrd}
\usepackage{url}

\newcommand\Lcal{\mathcal{L}}
\newcommand\Ncal{\mathcal{N}}
\newcommand\Bcal{\mathcal{B}}

\newcommand\Tcal{\mathcal{T}}
\graphicspath{{figures/}}
\usepackage{nomencl}
\usepackage[sort&compress]{natbib}
\usepackage{natbib}
 \bibpunct[, ]{(}{)}{,}{a}{}{,}%

 \providecommand{\keywords}[1]
{
  \small	
  \textbf{\textit{Keywords---}} #1
}

\floatname{algorithm}{Procedure}

\author[1]{Mertcan Yetkin\thanks{Corresponding Author: \texttt{mey316@lehigh.edu}}}
\author[1]{Brandon R. Augustino\thanks{\texttt{bra216@lehigh.edu}}}
\author[2]{Alberto J. Lamadrid\thanks{\texttt{ajlamadrid@lehigh.edu}, https://orcid.org/0000-0003-2539-9487}}
\author[1]{Lawrence V. Snyder\thanks{\texttt{lvs2@lehigh.edu}}}
\affil[1]{Department of Industrial and Systems Engineering, Lehigh University, Bethlehem, PA}
\affil[2]{Department of Economics, Lehigh University, Bethlehem, PA}

\date{\vspace{-5ex}}

\title{Co-optimizing the Smart Grid and Electric Public Transit Bus System}

\begin{document}

\maketitle

\begin{abstract}%
As climate change provides impetus for investing in smart cities, with electrified public transit systems, we consider electric public transportation buses in an urban area, which play a role in the power system operations in addition to their typical function of serving public transit demand. Our model considers a social planner, such that the transit authority and the operator of the electricity system co-optimize the power system to minimize the total operational cost of the grid, while satisfying additional transportation constraints on buses. We provide deterministic and stochastic formulations to co-optimize the system. Each stochastic formulation provides a different set of recourse actions to manage the variable renewable energy uncertainty: ramping up/down the conventional generators, or charging/discharging of the transit fleet. We demonstrate the capabilities of the model and the benefit obtained via a coordinated strategy. We compare the efficacies of these recourse actions to provide additional managerial insights. We analyze the effect of different pricing strategies on the co-optimization. We also conduct congestion analysis in the power network, comparing our cooperative approach to a non-cooperative strategy when we assume electrified fleet sizes grow with greater battery capacities. Given the recent momentum towards building smarter cities and electrifying transit systems, our results provide policy directions towards a sustainable future. We test our models using modified \textsc{MATPOWER} case files and verify our results with different sized power networks. This study is motivated by a project with a large transit authority in California. 
\end{abstract}


\section{Introduction}
\subsection{Motivation}

The global agenda established in the Paris agreement, which has been ratified by more than 180 countries  \citep{paris-agreement}, has declared climate change
to be at the forefront of challenges governments around the world are looking to overcome in the coming decades with green solutions. The global effort to reduce emissions encompasses various climate-conscious efforts, such as switching to renewable energy sources, decarbonization of transportation, and engineering of microbes to produce biofuels. Within the scope of electrification of mass transit systems, battery-electric buses play a crucial role in mitigation efforts. To this end, battery-electric buses (BEB) have been implemented throughout several countries including China, the United States, and many European countries \citep{pagliaro-electric-bus}. 

In this work, we address challenges related to the electrification of transportation systems associated with the joint operation of power and transit systems. The key question we address is to determine how a fleet of electric buses can best be integrated into public transportation infrastructure. In particular, the task at hand is to determine where and when to charge and discharge the buses such that operating is not a detriment to the power system, but rather, offers flexibility to the grid depending on the transportation--operational requirements. By co-optimizing, we ensure that the additional demand imposed by the transit authority is satisfied while inducing minimal stress on the power system.

\cite{pagliaro-electric-bus} argue that within the next decade, most fossil-fueled buses across the globe will be replaced by electric buses. \cite{Kelly934} discuss current transport solutions in place on air quality and emissions in London and Beijing. They conclude that emission-free vehicles are not enough to mitigate the problem, and efficient urban transportation can only be achieved by the expansion of mass transit systems. \cite{Creutzig911} discuss the effect of transportation on climate change mitigation. They note that many urban transport experts underline the reduction of overall transport demand growth via bus rapid transit and bicycle highways. They conclude that unless policy-makers' attitudes towards transport issues do not change soon, transport may hinder global efforts to mitigate climate change. \cite{chen-impacts-of-fleet-and-charging} study impacts of fleet types and charging strategies for electric vehicles on emissions in the Beijing area, where wind energy is considered to be the main generation source. They argue that utilizing electric buses and electric taxis is the most efficient way of curtailing NO$_\text{x}$ emissions. They underline that with the increase in wind penetration levels, if charging is not operated properly, electric vehicles could negate a significant amount of the environmental benefits obtained from renewables. Experts agree that emissions cannot be sufficiently reduced if industrialized countries continue their reliance on automobiles, even if those automobiles are fully electric. Instead, public transportation systems must be expanded and electrified. Our work contributes to this effort addressing potential challenges and benefits associated with day-to-day operation of electrified transit fleets. 

Along with the decarbonization of transit systems, proper coordination of different business and government branches is essential to ensure a feasible transition overall. \cite{Williams53} investigate infrastructure and technology paths required to meet California's reduction goals in emissions. They conclude that electrification of transportation, among other sectors, is essential to meet emission goals, and that infrastructure and technology must be deployed in a coordinated effort in order to realize their emissions-reduction goals. Further, the authors underline that the wide adoption of electric vehicles without smart charging (i.e. charging that is conscientious of the smart grid and the public) will impair dispatch operations. Their discussion motivates a system-wide coordinated effort within the charging problem which we address in this work.

 \cite{chu-opport-challenge} discuss opportunities and challenges for a sustainable energy future. They note that the integration of energy sources with electricity transmission, distribution, and storage is vitally important to offset the variations in renewable generation noting that ``Balancing and optimizing
power flow and generation are challenges that need improved technology, and business and regulatory vision." \cite{milovanoff-elec-light-duty-mitigation} investigate the sufficiency of the electrification of the light-duty fleet to reach mitigation targets in the US. They argue that cooperation between electric vehicle (EV) drivers and electricity suppliers is essential, and can be achieved by smart contracts. They highlight the need for policies to reduce vehicle ownership and usages such as the deployment of new public transit options, subsidies for public transportation, and transit-oriented land-use policies. In our work, we try to address some of the question put forth in policy-making. 

Coordination of operations handled by different parties not only increases the total efficiency of the system but also has the potential to introduce positive externalities such as dealing with cases when total generation exceeds demand (which potentially causes negative electricity prices).  \cite{zhou-negative-prices} analyze electricity trading strategies with potentially negative prices using a Markov decision process framework. They find that ignoring negative prices within the strategy can lead to a considerable loss of value. Further, they argue that under certain circumstances, the disposal of extra electricity purchased from the system with negative prices could be preferable compared to a storage strategy. Within the scope of our work, co-optimizing battery--electric buses and the smart grid would introduce a more efficient approach when electricity prices are negative. Specifically, the transit fleet could be co-optimized to store electricity from the system with negative prices, to later dispose of the electricity while they are serving transit demand. This additional flexibility would benefit both the ISO and the transit authority.

Building on the framework of the optimal power flow (OPF) problem  \citep{base-opf}, our models consider energy storage, uncertainty from renewable energy sources, vehicle-to-grid (V2G) capabilities, and the electrification of transit networks. We discuss the literature related to each of these aspects next.

\subsection{Operational aspects of co-optimizing the power and transit systems}
A subset of the topics studied in this paper have been extensively addressed in isolation in the literature (see e.g., \cite{smartTransit, chen2016optimal}). \cite{opf-management-with-pv-and-batteries} present an optimal power management scheme for photovoltaic (PV) systems with storage units. \cite{opf-with-storage} consider a microgrid with storage capabilities and aim to optimize the power consumption by balancing the power generation of renewable generators. However, these and other related papers differ from ours, as we aim to incorporate all of these aspects (renewables, storage, transit, etc.)\ in one formulation.

 The optimal power flow (OPF) problem, first studied by \cite{carpentier-opf}, is an optimization problem that determines the optimal dispatch in a power network, in which one solves for a network operating point that satisfies power flow equations and physical constraints such as thermal limits and line capacities  \citep{conejo2018power}. The multi-period optimal power flow problem (MPOPF) is a natural extension of the OPF problem, in which there are components involving multiple time periods such as storage units. In particular, storage units provide the capability of storing energy for a period of time and releasing it at a later time period, which introduces a temporal dependency rather than the snapshot of the system studied in the original OPF formulation. There are several works associated with the operation of the power grid by use of MPOPF, alternatively known as dynamic power flow (DOPF). By utilizing an MPOPF formulation, one can provide additional services to the power grid, such as voltage stability regulation and ramping reserves, which, in general, are referred to as {\em ancillary services} in the power systems literature; see, e.g.,  \cite{ancillary-services-opf}. The role of ancillary services in MPOPF is studied by \cite{mopf-vs-dr}, who propose an MPOPF formulation that incorporates demand-responsive loads specifically to improve steady-state voltage stability. Similarly, \cite{COSTA20071047} present a DOPF-based model to clear both energy and spinning reserve day-ahead markets. Moreover, \cite{lamadrid-ancillary} state that ancillary services are provided by uncertain renewable energy generation, which serves as ramping reserves. In the present work, we extend the standard MPOPF formulation to incorporate a fleet of public transportation buses that have the capability to charge and discharge within the power network.

\cite{mpopf-with-dr-and-uncerrenew} propose a stochastic MPOPF model to cope with the uncertainty stemming from renewable generators. A sparse formulation of a robust MPOPF problem with storage units and renewable generators is provided by \cite{robust-mpopf-with-storage-and-renew}, whose solution methodology utilizes receding horizon control. \cite{mpopf-with-wind} propose an MPOPF formulation that integrates wind farms into the generation portfolio. \cite{adaptive-robust-mpacopf} address the non-convexity in the power flow equations and uncertainty in renewable generation in terms of both active and reactive power. Their model incorporates transmission constraints and reactive capability curves of both conventional and renewable generators. To incorporate uncertain generation in our work, we also propose two additional formulations extended from stochastic MPOPF.

From a power systems perspective, our work can be considered as an extension of the MPOPF problem with (\textit{mobile}) storage units. Under this interpretation, the storage units are the buses in the transit fleet with additional operational constraints related to the transportation system. Many papers model the use of energy storage in the power grid. As a closely related work, a model similar to ours can be found in a study by \cite{mopf-with-storage}. Here, the authors consider an MPOPF formulation that incorporates large-scale standalone battery energy storage devices to reduce the fluctuations in the grid and to handle peak shaving in the grid. The objective in the formulation is to minimize the total generation cost (which can be different in each time period), including the charging and discharging costs of the batteries. They also argue that the batteries may be used to handle the uncertainty introduced when renewable energy sources are considered. Note, however, that no public transportation aspects are considered in their work. 

There are several papers in the transportation literature investigating problems associated with the electrification of transportation systems. \cite{charging_stockholm,charging_infra_model} propose mathematical optimization models to locate charging stations within an urban transportation network. They analyze the environmental outcomes, including emission of gases, based on a case study in Stockholm. 
\cite{energy_impact_auto_elec_fleet} investigate the effect of ambient temperature on the energy consumption of autonomous electric vehicles. Their results are demonstrated with a  data-driven simulated transportation network in New York City. Additionally, \cite{optimal-trafficpower-flow-in-urban-electrified-trans-network} focus on the coordinated operation of transportation and power systems. They assume an electrified transportation network capable of wireless power transfer coupled with a power network. They propose an optimal traffic--power flow model optimizing the generation schedule and congestion tolls as an optimization problem  
with traffic user equilibrium constraints.

The study of interdependent systems typically necessitates the use of multi-objective optimization techniques. \cite{moo-dist-energy-review} provide an overview of multi-objective planning of distributed energy resources (DERs) where the objective function involves terms from perspectives such as the distributed energy resources developer, the distribution system operator, and the regulator. 
They underline the importance of DER integration and argue that poor integration of DERs can result in increases in losses, as well as voltage and network instability. The insights provided by \cite{moo-dist-energy-review} exemplify the benefits obtained by a co-optimization framework. In our case, the co-optimization is considered between the independent system operator (ISO) of the power network and the transit authority.

From the perspective of interdependent systems, only a few papers focus on the joint operation of electric vehicles and the power grid. \cite{mult-obj-dopf-ev} aim to coordinate plug-in electric vehicles in a multi-objective security-constrained DOPF problem to minimize the total operation cost and emissions.  \cite{mult-obj-ev-charge} propose a multi-objective method for charging/discharging electric vehicles to minimize the total operational costs and emissions. 
An important assumption in their methodology is that the EV owner decides the parking time for charging/discharging one day in advance. This assumption is one of the key differences between public and private means of transportation. 
Because it considers public transportation systems, our study assumes that the decision of when and where to charge/discharge is made as a consequence of the joint operation between the independent system operator (ISO) and the public transportation authority (TA).

A related methodology to ours is by \cite{lin-shen-charging-stations}, who study a planning problem that decides the locations and sizes of charging stations, considering a coupled transportation network and power network. Specifically, they seek to determine a long-term plan with a horizon of roughly 10 years. An immediate drawback is that they assume that all information (electricity prices, infrastructure costs, demand, technology) will be relatively unchanged over their planning horizon. Thus, having such a long-term strategy may not be preferable. Further, the OPF problem is typically solved many times a day, sometimes as often as every five minutes \citep{cain2012history}, whereas the formulation by \cite{lin-shen-charging-stations} couples the transportation aspects with OPF only at two periods. The first of these represents the first 10 years of the planning horizon while the second stage is the once at a given stage of length 10 years, and the second represents 30 years. Another shortcoming of their model is that it captures instantaneous {\em power} flow, whereas storage optimization requires calculations of {\em energy} (power accumulated over time). In contrast, the time dimension in our model (and other MPOPF-based models) captures the relationship between power and energy. 

From the perspective of public transit operators, \cite{abdelwahed-charge-transit-elec-bus} develop models 
to optimize the charging operation of a fleet in Rotterdam. However, they do not address the coupling effect of the power system, such as power system operation and limitations, V2G possibilities, and offset of the intermittency of renewable generation.

In sum, our work differs from the related works in the literature in the following aspects: we consider \textit{(i)} a fleet of electric vehicles and a power grid operated by a social planner, \textit{(ii)} operational constraints related to the transit fleet while providing services to the power grid, including V2G, \textit{(iii)} schedules for the battery (transit bus) connection, and \textit{(iv)} relocation of the batteries within the power grid. Further, we address the joint operation of an electric public transit system and the power system, specifically when the fleet is ``off-schedule. We use ``off-schedule" to refer to the set of time periods when the buses are not on their routes.


Our primary contributions are as follows:
\begin{enumerate}
    \item To the best of our knowledge, this is the first paper to address this problem in detail. We note that a preliminary, and abbreviated, version of this work appeared in \cite{Yetkin-Augustino-abstract}.
    \item We provide a deterministic formulation that jointly optimizes the operation of a public transit authority and an ISO.
    \item In the presence of renewable generation, we further extend this model to two two-stage stochastic programs (2SSPs) with different recourse actions, including additional charging/discharging of the transit fleet and ramping up/down the conventional generators.
    \item We additionally provide managerial insights: first, by conducting a benefit analysis of coordinated optimization; second, by demonstrating potential benefits of alternative recourse actions in the presence of variable renewable energy uncertainty; third, by analyzing the effect of pricing on the outcomes of co-optimization; and fourth, by conducting a congestion analysis when we consider growing fleet battery sizes.
\end{enumerate}

The rest of the paper is organized as follows. In Section~\ref{sec:problem_state} we provide a formal definition of the problem. In Section~\ref{sec:model_det} we introduce the deterministic formulation of the optimization problem. In Sections~\ref{sec:stoc_model_ramp} and \ref{sec:stoc_model_charge} the two 2SSP formulations are provided. Numerical results are presented in Section~\ref{sec:num}, managerial insights are provided in Section~\ref{sec:managerial_insights}, and Section~\ref{sec:con} concludes the paper. Additional content can be found in Appendices~\ref{apx:nomenclature}--\ref{apx:MPOPF}. 

\section{Problem statement}
\label{sec:problem_state}
We consider a single social planner who manages both the power and public transit systems. The goal of the social planner is to co-optimize the joint operation of these two systems, by optimizing the charging/discharging of the transit bus batteries when buses are off-schedule, over a horizon of one day. The only operational requirement on the transit buses is that they have to be fully charged before starting their schedules on the following day. The following decisions are addressed in the formulation while satisfying operational constraints: \textit{(i)} where, among a prespecified subset of nodes that serve as charging stations, to locate transit buses to charge/discharge, \textit{(ii)} how much electricity to charge/discharge, and \textit{(iii)}  power dispatch.

We simplify the problem by assuming the following:
\begin{itemize}
    \item All of the information related to the power network is known (topology, generator limits, line limits, etc.), except the ramping costs of conventional generators (which will be investigated in Section~\ref{sec:benefit_stoc}). \item Over the course of one day (the planning horizon) the topologies of both the transportation and power networks remain fixed. The power network is provided as a graph with edges representing power lines and nodes representing connection points of lines. In the transportation network, the nodes of the graph are interpreted as charging stations, and the edges as the transit routes among charging stations. Note that the charging stations in the transportation network act as coupling points between the two networks since they also appear as a subset of nodes in the power network.
    \item Electricity demand is known for the entire horizon, with the further assumption that an accurate day-ahead estimation of the demand is accessible.
    \item We consider direct current (DC) approximation for the power system  \citep{conejo2018power}. This assumption is useful for simplifying the computational complexity of our optimization problem since the alternating current (AC) OPF problem is known to be nonlinear and non-convex \citep{GOPINATH2020106688}.
    \item We assume black start capabilities for conventional generators. This assumption eliminates the start-up cost of conventional generators.
    \item All of the information related to the transit system is known (i.e. transportation network, schedules, travel times).
    \item All of the information related to batteries of the transit vehicles are known (the charge/discharge rate, total capacity, efficiency). We assume that the efficiency of the battery is the same for charging and discharging, though this assumption can be easily relaxed. 
    \item Each charging station in the network has enough capacity for the entire fleet to charge simultaneously at one location. A related discussion is also provided in Appendix~\ref{apx:det_model}. Moreover, buses can discharge at a charging station as part of the V2G technology.
\end{itemize}

\section{Deterministic formulation}  \label{sec:model_det}
In this section, we present the deterministic formulation to co-optimize the operation of the transit fleet and the power grid. The formulation can be seen as an extension of the MPOPF problem with additional transportation aspects. The problem is a mixed-integer quadratic program (MIQP). The following components are captured by this formulation:
\begin{itemize}
    \item Ensuring each bus has a full battery at the end of their off-schedule (also considered by \cite{abdelwahed-charge-transit-elec-bus}).
    \item Dispatch of power subject to physical constraints of the grid, such as ramping rates, transmission limits, etc.
    \item Charging/discharging of the buses subject to physical constraints with respect to the batteries of the transit buses, such as charging capacity and charging/discharging rate.
    \item Relocation of the buses subject to constraints related to the travel time of the buses depending on the node to connect to the electrical grid.
\end{itemize}
Due to space limitations, we provide the complete deterministic formulation in Appendix~\ref{apx:det_model}, and we verbally describe the model in this section. 

Our objective aims to minimize a convex combination, with coefficient $\alpha$, of the total power generation cost and the charging/discharging cost of the transit buses. 
The model contains standard DC OPF constraints with adjustments to include the effects of charging and discharging. Then, based on the fleet's schedule, we ensure that initial and final battery levels are respected. The model also includes updates for the battery levels. Assignment constraints ensure that a bus must either be traversing the network or be stationed at a connection point for charging/discharging at any time point. Specifically, a bus can only relocate between two points in the power network if there is enough time available. We initially assign all of the fleet to a depot node. In addition, we have several bounds on ramping amounts of generators, generation amounts, line flows, charge/discharge amounts of batteries, and battery levels.

The two terms in the objective function each ``belong'' to a different party: The generation costs are incurred by the power grid operator while the charging costs are incurred by the transit operator. We use a convex combination of these two terms since they are usually in different scales, and since the central decision-maker may wish to place more or less weight on one term or the other. Similarly, in the stochastic models presented in the following sections, we use a convex combination of costs that are in different scales, to differentiate between the costs of different parties in this cooperation. In all of our numerical experiments, we set the convex combination coefficient, $\alpha$, to 0.5.
Based on the requirements of the specific application, one can examine different objective terms such as line losses, voltage deviations, and many others tailored towards the desired goal. Note that, due to the specific application, we consider off-schedule periods as a single block of time (e.g., 5 pm--4 am); if there are multiple, non-contiguous blocks, one should account for initial and final conditions on the battery levels for each block of off-schedule times. 

\section{Two-stage stochastic formulations}
In this section, we consider renewable generation units (e.g., wind generators) within the power grid. In each time period, the quantity of renewable generation is random. We incorporate the uncertainty regarding the renewable generation units using scenarios in a two-stage stochastic modeling approach. Each of these scenarios represents one realization of wind generation distribution. The details of the scenario generation scheme are described in Section~\ref{sec:results_stoch}.
In addition to the assumptions listed in Section~\ref{sec:problem_state}, we assume that a forecast of the renewable generation is available for one day, though the actual generation may differ randomly from the forecast.

Specifically, we present two different two-stage stochastic formulations, which assume different recourse actions. In the first formulation (Section~\ref{sec:stoc_model_ramp}), the recourse action is ramping up/down the conventional generators, whereas in the second formulation (in Section~\ref{sec:stoc_model_charge}), the recourse is additional charging/discharging of the transit fleet. The first formulation assumes that the ISO is taking the recourse actions to mitigate the uncertainty, whereas the second assumes the transportation authority is doing so.

\subsection{Ramping-based formulation}
\label{sec:stoc_model_ramp}
The following formulation can be considered as an extension of the two-stage stochastic single-period OPF formulation presented by \cite{morales2013integrating} in the sense that we have the additional aspect of transportation and related operational constraints in the first stage. Briefly, we have the following set of decisions in addition to those listed in Section~\ref{sec:model_det}: a first-stage commitment of the renewable generation units on how much to generate; and a second-stage adjustment of power generation via ramping up/down the conventional generators. We again omit the complete model and provide the details in Appendix~\ref{apx:stoch_ramping_model}.

Our objective is the summation of first-stage costs (as in the deterministic objective function) and second-stage costs including the expected renewable generation costs and expected ramping costs. Added to the constraints captured in the deterministic model in Section~\ref{sec:model_det} are the second-stage nodal balance and flow constraints. We are further constrained by limits on renewable generation, limits on second-stage ramping, and limits on shedding in the second stage. We also have constraints on the transit bus battery levels, charging limits, and transit fleet operation. More importantly, we have each of the OPF constraints repeated for multiple time periods rather than a single time period, in contrast to \cite{morales2013integrating}, where the time periods in our formulation are coupled via batteries on the transit fleet. The inclusion of multiple periods in the formulation increases the complexity of the problem dramatically.

\subsection{Charging/discharging-based formulation}
\label{sec:stoc_model_charge}

Rather than handling the uncertainty by ramping up/down the conventional generators, we now consider this service to be handled by the transportation authority via additional charging/discharging of the transit vehicles in the second stage. We assume that the realization of the scenarios occurs in near-real-time so that there is not enough time to relocate the transit buses after the scenario realizations. Then, we have the following set of decisions in addition to the ones in Section~\ref{sec:model_det}: a first-stage commitment of the renewable generation units on how much to generate; and an additional second-stage charging/discharging of the transit fleet.

A complete table of nomenclature can be found in Appendix~\ref{apx:nomenclature}. The formulation of the 2SSP whose recourse is charging/discharging of the transit buses is given as follows:
\begin{align}
    \min \quad &(1-\alpha) \left[ \sum_{t \in \Tcal} \sum_{i \in \Ncal_g} c^g_{it} p^g_{i t} + c'^g_{it} (p^g_{i t})^2 + \sum_{\omega \in \Omega} \pi_\omega \left ( \sum_{t \in \Tcal} \sum_{i \in \Ncal_r} c^r_{it} p^r_{it\omega}  +  \sum_{t \in \Tcal}\sum_{i \in \Ncal }c_{it}^{\text{shed}} p_{it\omega}^{d,\text{shed}} \right )  \right ]  \nonumber \\
    & + \alpha \left [ \sum_{t \in \Tcal} \sum_{b \in \Bcal} \sum_{i \in \Ncal_b} c_{i t} \left(p^c_{ibt} - p^{dc}_{ibt} \right) + \sum_{\omega \in \Omega }\pi_\omega  \sum_{t \in \Tcal} \sum_{b \in \Bcal} \sum_{i \in \Ncal_b} c_{i t \omega}^{+} \left(p^{c,+}_{ibt\omega} - p^{dc,+}_{ibt\omega} \right )  \right ] \label{eq:obj1_s2}
\end{align}
subject to:
\begin{alignat}{2}
    & p^g_{it} + p^r_{it} - \sum_{b \in \Bcal} p^c_{ibt} + \sum_{b \in \Bcal} p^{dc}_{ibt} - p^d_{it} = \sum_{j : (i,j) \in \mathcal{L}} p_{ijt} - \sum_{j : (j,i) \in \mathcal{L}} p_{jit}
    && ~\forall i \in \Ncal, t \in \Tcal \label{eq:ssc1} \\ 
    & p^r_{it\omega}-p^r_{it} - \sum_{b \in \Bcal} p^{c,+}_{ibt\omega} + \sum_{b \in \Bcal} p^{dc,+}_{ibt\omega} + p^{d,shed}_{it\omega}= \sum_{j : (i,j) \in \mathcal{L}} \left (p_{ijt\omega} - p_{ijt} \right ) - \sum_{j : (j,i) \in \mathcal{L}} \left ( p_{jit\omega} - p_{jit}\right )
    && ~\forall i \in \Ncal, t \in \Tcal, \omega \in \Omega \label{eq:ssc2}
    \end{alignat}
    \begin{align}
     &p^g_{it} = 0  &\forall i \in \Ncal \setminus \Ncal_g, t \in \Tcal \label{eq:ssc3} \\
     &p^r_{it} = 0  &\forall i \in \Ncal \setminus \Ncal_r, t \in \Tcal \label{eq:ssc4} \\
     &p^r_{it\omega} = 0  &\forall i \in \Ncal \setminus \Ncal_r, t \in \Tcal, \omega \in \Omega \label{eq:ssc5} \\
     & -\overline{\theta} \leq \theta_{it} \leq \overline{\theta} & \forall i \in \Ncal, t \in \Tcal \label{eq:ssc6} \\
     & -\overline{\theta} \leq \theta_{itw} \leq \overline{\theta} & \forall i \in \Ncal, t \in \Tcal, \omega \in \Omega \label{eq:ssc7} \\
     & \theta_{1t} = 0 &\forall t \in \Tcal \label{eq:ssc8} \\
     & \theta_{1t\omega} = 0 &\forall t \in \Tcal, \omega \in \Omega \label{eq:ssc9} \\
     & p_{ijt} = \frac{\theta_{it} - \theta_{jt}}{x_{ij}} & \forall (i,j) \in \Lcal, t \in \Tcal \label{eq:ssc10} \\
       & p_{ijt\omega} = \frac{\theta_{it\omega} - \theta_{jt\omega}}{x_{ij}} & \forall (i,j) \in \Lcal, t \in \Tcal, \omega \in \Omega \label{eq:ssc11}\\
    & -\overline{S}_{ij} \leq p_{ijt} \leq \overline{S}_{ij} & \forall (i,j) \in \Lcal, t \in \Tcal \label{eq:ssc12}\\
    & -\overline{S}_{ij} \leq p_{ijt\omega} \leq \overline{S}_{ij} & \forall (i,j) \in \Lcal, t \in \Tcal, \omega \in \Omega \label{eq:ssc13}\\
    &0 \leq p^g_{it} \leq \overline{p}_{it}^g &\forall i \in \Ncal_g, t \in \Tcal \label{eq:ssc14}\\
    & -p^{\delta g}_{it} \leq p^g_{i,t+1} - p^g_{it}\leq p^{\delta g}_{it} &\forall i \in \Ncal_g, t \in \Tcal \setminus \{T\} \label{eq:ssc15}\\
    &0 \leq p^r_{it} \leq \overline{p}_{it}^r &\forall i \in \Ncal_r, t \in \Tcal \label{eq:ssc16}\\
    &0 \leq p^r_{it\omega} \leq \tilde{p}_{it\omega}^r &\forall i \in \Ncal_r, t \in \Tcal, \omega \in \Omega \label{eq:ssc17}\\
    &0 \leq p^{d,\text{shed}}_{it \omega} \leq p^d_{it} & \forall i \in \Ncal, t \in \Tcal, \omega \in \Omega \label{eq:ssc18}\\
     &e_{b T^1_b \omega} = e^1_b &\forall b \in \Bcal, \omega \in \Omega \label{eq:ssc19}
    \end{align}
    \begin{equation}
        e_{b T^2_b \omega} + \eta_b \sum_{i \in \Ncal_b} \left ( p^c_{ib T^2_b} + p^{c,+}_{ib T^2_b \omega } \right ) \delta t - \frac{1}{\eta_b} \sum_{i \in \Ncal_b} \left ( p^{dc}_{ibT^2_b} + p^{dc,+}_{ibT^2_b \omega} \right )\delta t - s_b  y_{bT^2_b} = \overline{e}_b ~~\hspace{1cm} \forall b \in \Bcal, \omega \in \Omega \label{eq:ssc20}
    \end{equation}
    \vspace{-0.5cm}
    \begin{multline}
        e_{b,t+1,\omega} = e_{bt \omega} + \eta_b \sum_{i \in \Ncal_b} \left ( p^c_{ibt} + p^{c,+}_{ibt\omega} \right )\delta t - \frac{1}{\eta_b} \sum_{i \in \Ncal_b} \left ( p^{dc}_{ibt} + p^{dc,+}_{ibt\omega}\right) \delta t - s_b  y_{bt} \\~~ \forall b \in \Bcal, \: t,t+1 \in \Tcal_b, \omega \in \Omega \label{eq:ssc21}
    \end{multline}
    \begin{align}
    &\underline{e}_{b} \leq e_{bt\omega} \leq \overline{e}_{b}  &\forall b \in \Bcal, t \in \Tcal_b, \omega \in \Omega \label{eq:ssc22}\\
    & p^c_{ibt} + p^{c,+}_{ibt\omega} \leq \overline{p}^c_{b} z_{ibt} &\forall i \in \Ncal_b, b \in \Bcal, t \in \Tcal, \omega \in \Omega \label{eq:ssc23}\\
    &p^{dc}_{ibt} + p^{dc,+}_{ibt\omega} \leq \overline{p}^{dc}_{b} z_{ibt} &\forall i \in \Ncal_b, b \in \Bcal, t \in \Tcal, \omega \in \Omega \label{eq:ssc24}\\
    &p^c_{ibt} \geq 0, ~ p^{dc}_{ibt} \geq 0 &\forall i \in \Ncal_b, b \in \Bcal, t \in \Tcal \label{eq:ssc25}\\
    &p^{c,+}_{ibt\omega} \geq 0, ~ p^{dc,+}_{ibt\omega} \geq 0  &\forall i \in \Ncal_b, b \in \Bcal, t \in \Tcal, \omega \in \Omega \label{eq:ssc26}\\
     &p^c_{ibt} = 0, ~ p^{dc}_{ibt} = 0 &\forall i \in \Ncal \setminus \Ncal_b, b \in \Bcal, t \in \Tcal \label{eq:ssc27}\\
     &p^{c,+}_{ibt\omega} = 0, ~ p^{dc,+}_{ibt\omega} = 0 &\forall i \in \Ncal \setminus \Ncal_b, b \in \Bcal, t \in \Tcal, \omega \in \Omega \label{eq:ssc28}\\
     &\sum_{i \in \Ncal_b} z_{ibt} + y_{bt} = 1&\forall b \in \Bcal, t \in \Tcal_b \label{eq:ssc29}\\
     &z_{ibt} + z_{jbt'} \leq 1 &  \forall t' \in \Tcal_b, ~ t < t' \leq t+\Delta t(i,j), \nonumber \\
     & & \forall i,j \in \Ncal_b, i \neq j,  b \in \Bcal, t \in \Tcal_b \label{eq:ssc30}\\
     & z_{d b T^1_b} = 1 &  d = \Ncal_{b}(1), \forall b \in \Bcal \label{eq:ssc31}\\
     &y_{bt} = 0 &\forall b \in \Bcal, t \in \Tcal \setminus \Tcal_b \label{eq:ssc32}\\
     &z_{ibt} = 0 &\forall i \in \Ncal_b, b \in \Bcal, t \in \Tcal \setminus \Tcal_b \label{eq:ssc33}\\
     &y_{bt} \in \{0,1\} & \forall b \in \Bcal, t \in \Tcal \label{eq:ssc34}\\
     &z_{ibt} \in \{0,1\} & \forall i \in \Ncal_b, b \in \Bcal, t \in \Tcal \label{eq:ssc35}
\end{align}

The objective function \eqref{eq:obj1_s2} accounts for the conventional generation cost in the first stage, expected renewable generation cost, charging/discharging cost of the fleet in the first stage and expected charging/discharging cost of the fleet in the second stage. We have already introduced the following sets of constraints in the deterministic formulation in Section~\ref{sec:model_det}: first-stage DC optimal power flow constraints \eqref{eq:ssc1}, \eqref{eq:ssc3}, \eqref{eq:ssc4}, \eqref{eq:ssc6}, \eqref{eq:ssc8}, \eqref{eq:ssc10}, \eqref{eq:ssc12}; bounds on generation \eqref{eq:ssc14}, ramping \eqref{eq:ssc15}, and first-stage charging/discharging \eqref{eq:ssc25}, \eqref{eq:ssc27}; assignment constraints for vehicle to charging station connection and vehicle relocation \eqref{eq:ssc29}--\eqref{eq:ssc35}.

We now additionally have second-stage nodal balance and flow constraints \eqref{eq:ssc2},  \eqref{eq:ssc7},  \eqref{eq:ssc9},  \eqref{eq:ssc11}, \eqref{eq:ssc13}; limits on renewable generation \eqref{eq:ssc16}, \eqref{eq:ssc17}; and limits on shedding in the second stage \eqref{eq:ssc18}. Unlike the formulation discussed in Section~\ref{sec:stoc_model_ramp}, we are also charging/discharging in the second stage. Thus, the battery level variables for vehicles $e_{bt\omega}$ are second-stage variables over the set of scenarios $\omega \in \Omega$. Also, update constraints for battery levels are replicated for each scenario in \eqref{eq:ssc19}--\eqref{eq:ssc22}. Moreover, we have modified bounds on the charging/discharging of the vehicles in \eqref{eq:ssc23}, \eqref{eq:ssc24}. Note that the positioning variables for the vehicles $z_{ibt}, y_{bt}$ remain as first-stage variables, since we assume that there is insufficient time to adjust the location of the vehicles after the scenario realizations.

\section{Numerical experiments}\label{sec:num}

In this section, we first provide (in Section~\ref{sec:results_det}) the results of a small case study for the deterministic formulation, then (in Section~\ref{sec:results_stoch}) the results of the stochastic formulations. Based on the transportation--electrification study on which we have collaborated with Santa Clara Valley Transportation Authority (\href{https://www.vta.org}{VTA}), we will focus on a regional study around San Jose, CA. Our partner VTA has been undertaking the integration of electric transit buses into their fleet in the city of San Jose, initially for a selected number of routes within their operating region.

\subsection{Case study for the deterministic formulation}

\label{sec:results_det}

\begin{figure}
\begin{center}
\begin{subfigure}{.69\textwidth}
\includegraphics[trim=0.5cm 0.5cm 1.5cm 0.5cm,width=0.6\textwidth,angle =270]{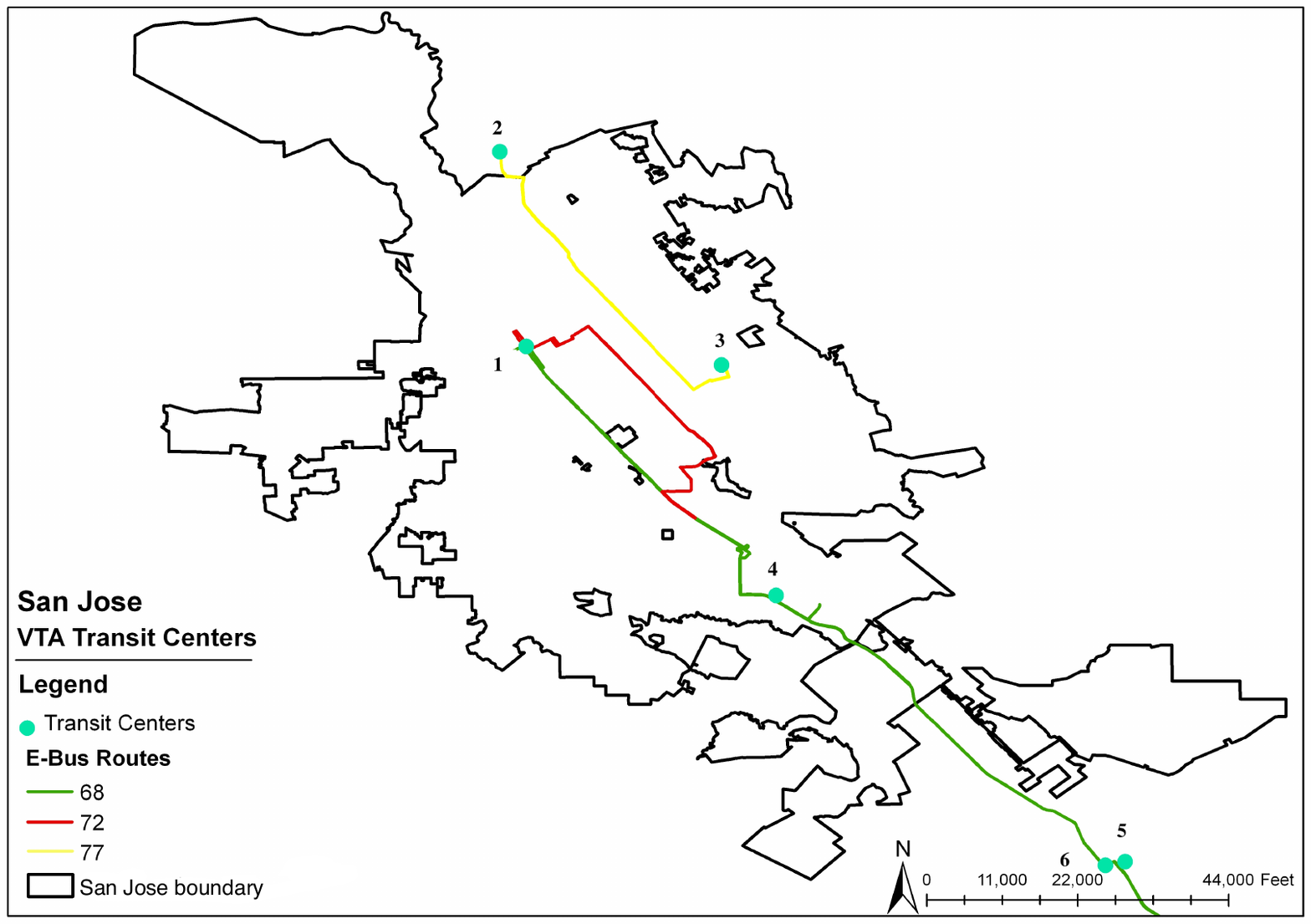}
\end{subfigure}%
\begin{subfigure}{.29\textwidth}
\includegraphics[trim= 7cm 5cm 5cm 5cm,scale=0.3]{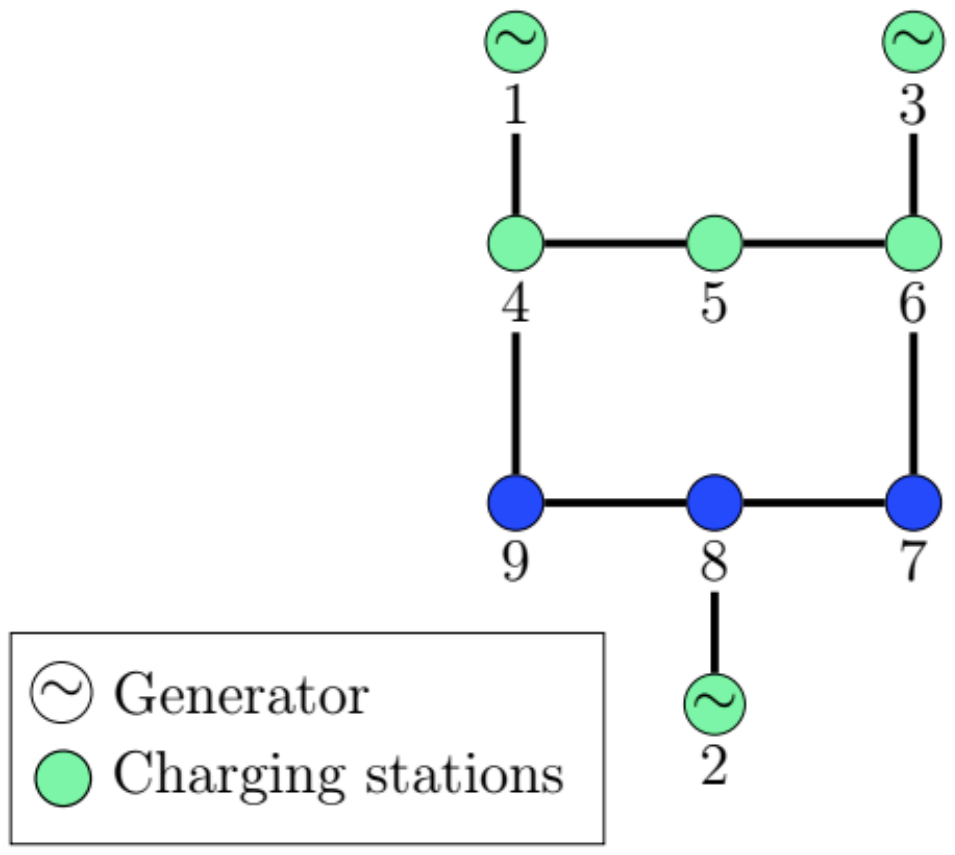}
\end{subfigure}
\end{center}
\caption{Transit layout (left) and power network schematic (right), where charging stations are assumed to be located at Transit Centers.}
\label{fig:transit}
\end{figure}

We consider a case study consisting of synthetic but realistic data meant to reflect the actual transit network in San Jose, CA. We overlay the 9-bus power network from \textsc{MATPOWER} \citep{MATPOWER} atop the geographical area.  Related parameters including demand ($p_{it}^d$), line limits ($\overline{S}_{ij}$), generation limits ($\overline{p}^g_{it}$) and cost of generation ($c^g_{it},c'^g_{it}$) are obtained from the \textsc{MATPOWER} case file. Figure \ref{fig:transit} provides a visualization of the layout, indicating the locations of charging stations for electric bus connection in both the power and transit networks. Note that the transportation network is smaller than the power network since the charging stations are chosen to be a subset of nodes in the power network. In general, a power network is not confined within a city limit and hence is likely to encompass the transit region.

\begin{figure}
\begin{center}
\includegraphics[width=8cm]{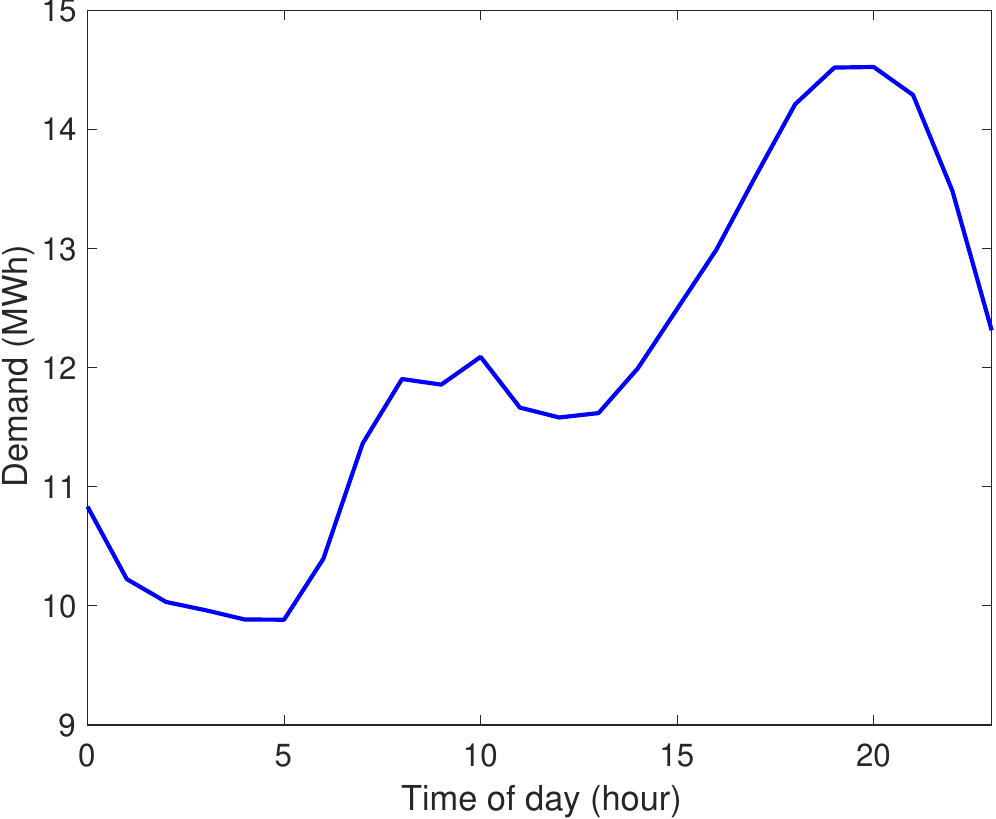}
\end{center}
\caption{Daily demand data}
\label{fig:demand_data}
\end{figure}

\begin{figure}
\begin{center}
\includegraphics[width=10cm]{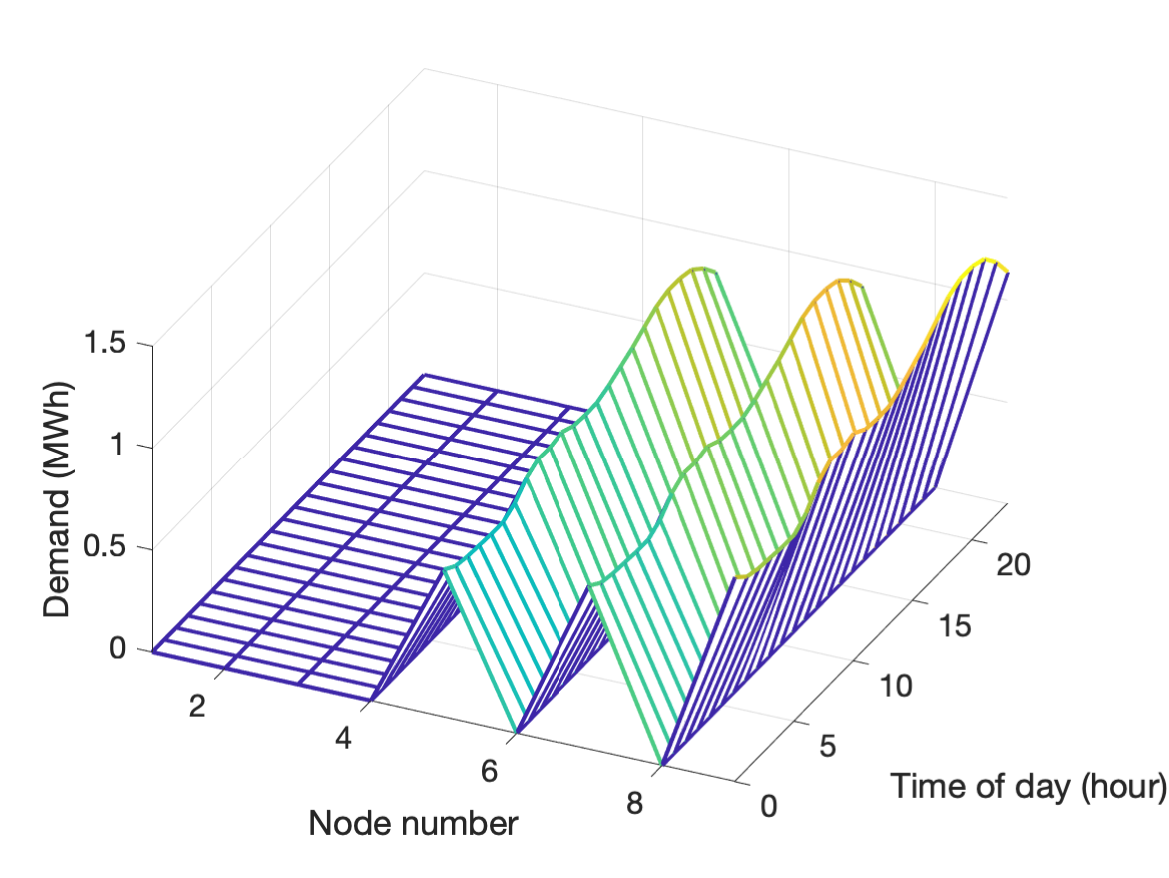}
\end{center}
\caption{Complete demand profile in power network (constructed)}
\label{fig:complete_demand_prof}
\end{figure}

\begin{figure}
\begin{center}
\includegraphics[width=10cm]{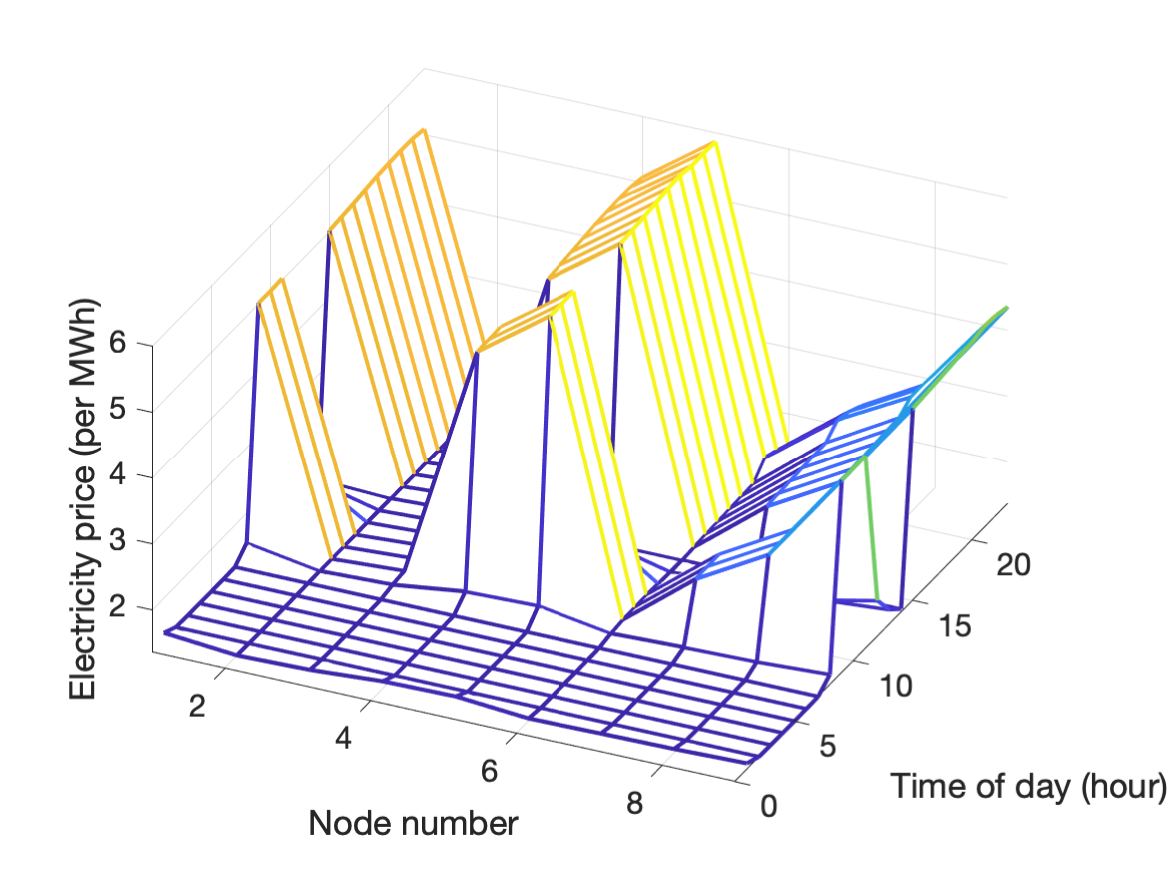}
\end{center}
\caption{Price data}
\label{fig:price_data}
\end{figure}

Line limits and demands from the 9-bus system are scaled down to the order of around 1 MWh so that we can clearly analyze the impact of the electric buses on the power network since battery energy capacities considered in this work are 0.66 MWh. Voltage angle limits $\overline{\delta}$ are set to $\pi/2$ and ramping limits are chosen to be $\frac{\overline{p}^g_{it}}{5}$.
Since the standard case file in \textsc{MATPOWER} provides a snapshot of the system in one time period, we expand the demand profile over the multiple time periods based on information obtained from the California Independent System Operator (\href{https://www.caiso.com}{CAISO}) for the San Jose region.  Moreover, the charging/discharging costs ($c_{it}$) in the objective function are the electricity prices obtained by averaging the prices as described in Section~\ref{sec:benefit_det}. The total daily demand and the price data are displayed in Figures \ref{fig:demand_data} and \ref{fig:price_data}, respectively. We also provide the expanded demand profile for the power network over multiple time periods in Figure~\ref{fig:complete_demand_prof}.

We assume that the four transit buses are 40-foot Proterra Catalyst E2Max models; the data on battery consumption, battery capacity, and charge/discharge limits are obtained from \href{https://www.proterra.com}{Proterra}. The battery efficiency is set to 0.9.
We consider hourly time steps and a 24-hour horizon. The schedules of the buses and charging station locations are obtained from VTA, where $\Tcal_b$ represents the time periods in which transit bus $b$ is in its off-schedule. $\Delta t(i,j)$ values are travel times between charging stations calculated using Google Maps (and discretized based on the resolution of the formulation).

We performed the optimization using Gurobi version 8.1.0 developed by \cite{gurobi}, with the default settings. The hardware was 2.2 GHz Dual-Core Intel Core i7 and 8GB memory. The solution time for this model was approximately 23 seconds.

\begin{figure}
\begin{center}
\begin{subfigure}{.19\textwidth}
  \includegraphics[trim = 0cm 0cm 0cm -0.5cm, scale=.35]{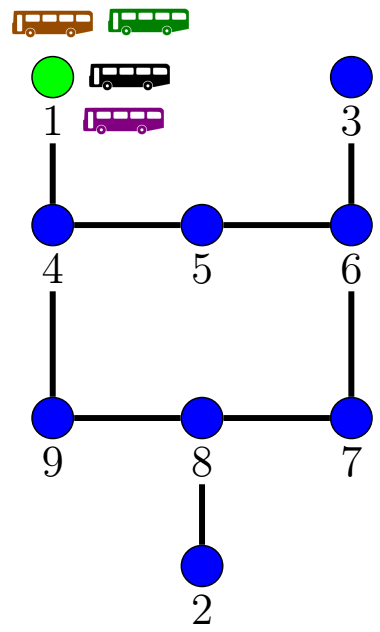}
\end{subfigure}
\begin{subfigure}{.19\textwidth}
  \includegraphics[trim = 0cm 0cm 0cm -0.5cm,scale=.35]{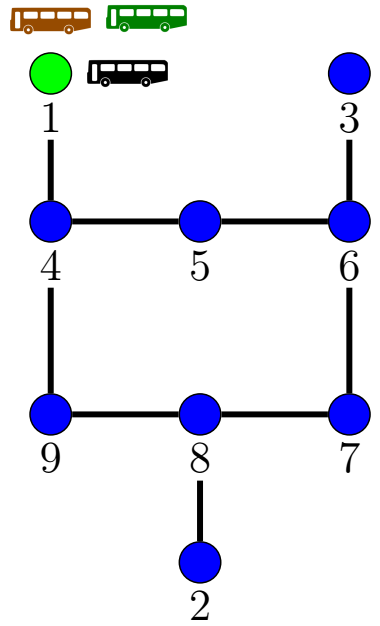}
\end{subfigure}
\begin{subfigure}{.19\textwidth}
  \includegraphics[trim = 0cm 0cm 0cm -0.5cm,scale=.35]{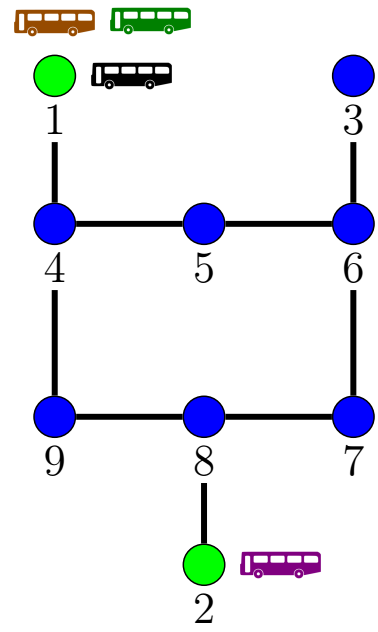}
\end{subfigure}
\end{center}
\caption{Transit bus locations for $t = 21$, $t = 22$ and $t = 23$ in optimal solution}
\label{fig:locations}
\end{figure}

\begin{figure}
\begin{center}
\begin{subfigure}{.43\textwidth}
  \includegraphics[width=7cm]{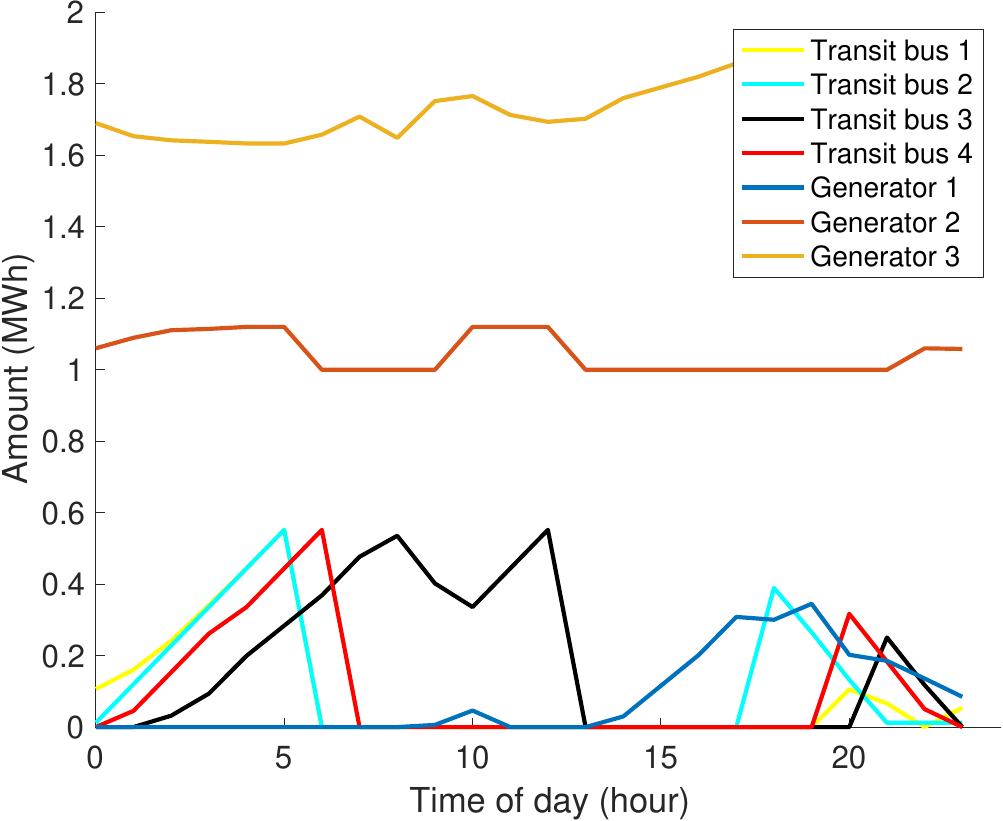}
\end{subfigure}
\begin{subfigure}{.43\textwidth}
    \includegraphics[width=7cm]{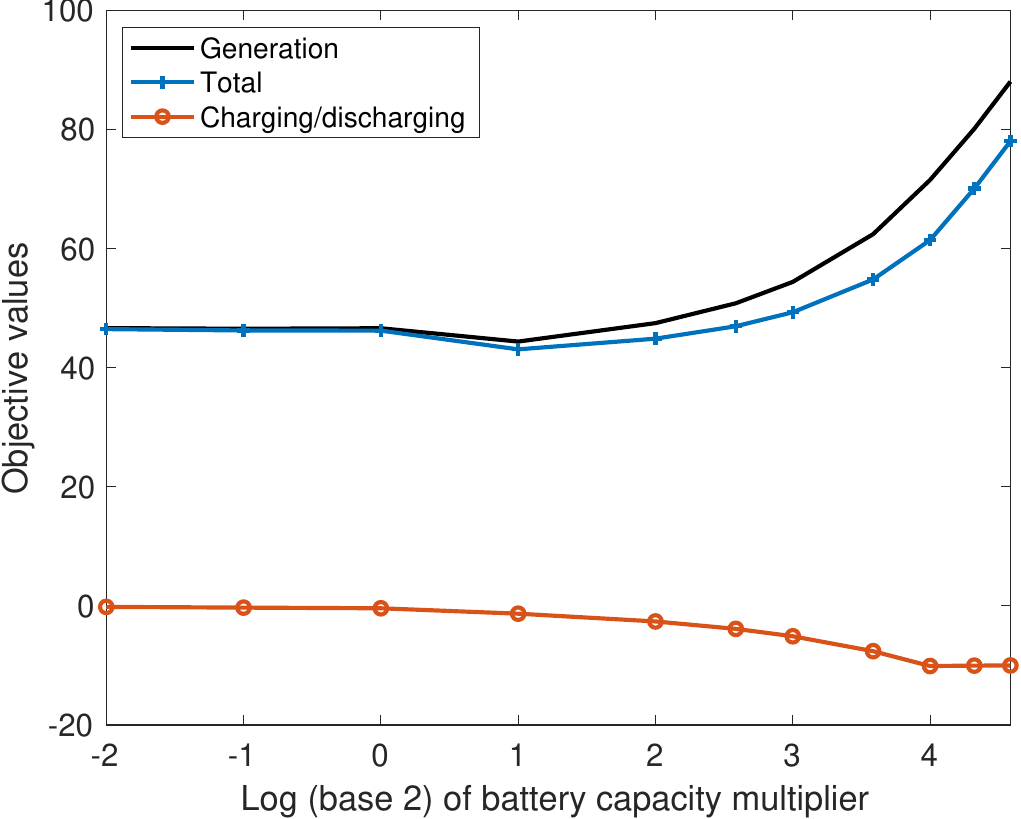}
\end{subfigure}
\end{center}
\caption{Generation and battery level profile (left) and objective values via varying battery capacity (right)}
\label{fig:solution}
\end{figure}

Figure~\ref{fig:locations} displays the optimal locations of the transit buses for 3 consecutive time periods. At time $t=21$, all of the four vehicles are located at node 1, the transit depot. In the next time period, three of the buses remain at node 1, while the fourth is in transit. 
At $t=23$, this bus arrives at its new location, node 2, to decrease total generation cost.
Figure~\ref{fig:solution} (left) shows the generation and battery level profile. In the optimal solution, battery levels in the first periods increase due to the operational constraints and decrease in the last periods since the price of electricity and the demand are high. Further, during the periods in which the transit buses are on their schedule, we can observe that the battery levels (i.e., the amount of energy stored in the batteries of the transit vehicles when they are not connected to the grid) are displayed as zero.

\begin{figure}
\begin{center}
\includegraphics[width=10cm]{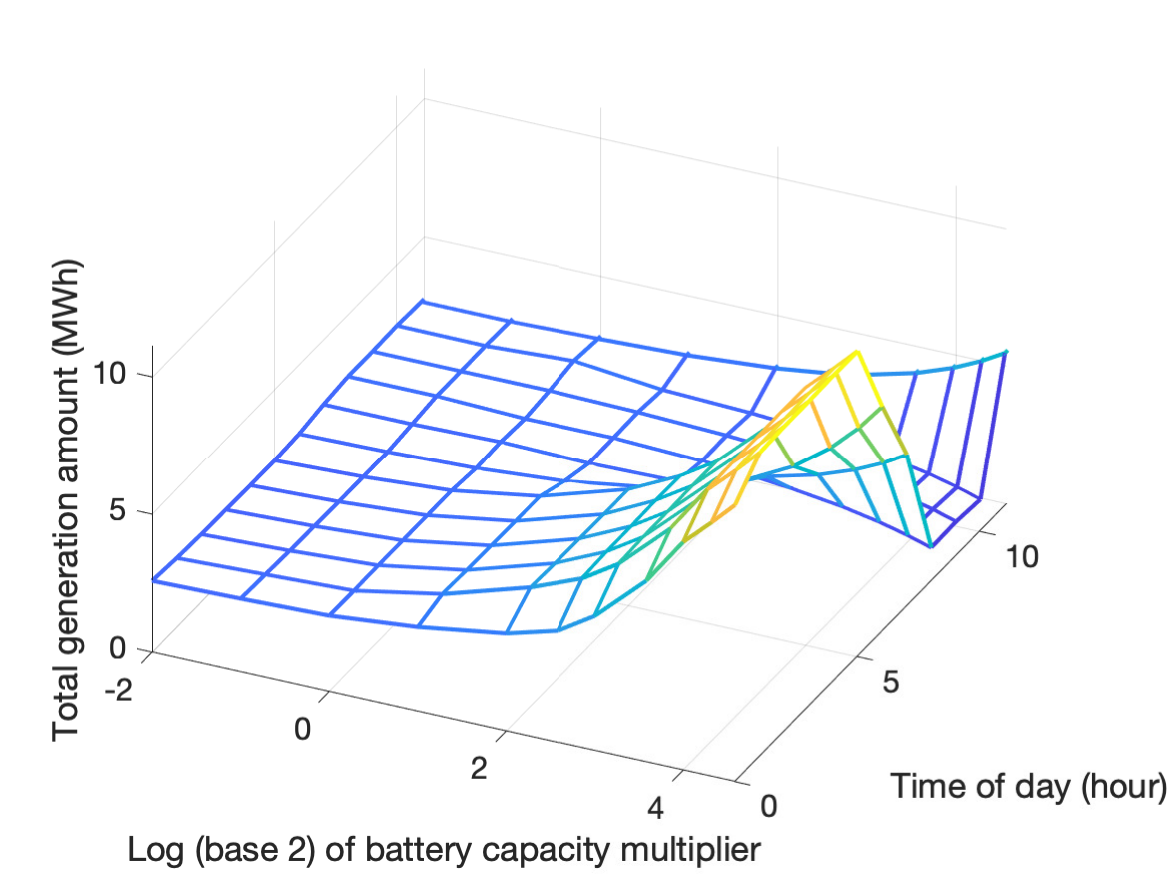}
\end{center}
\caption{Total generation via varying battery capacity}
\label{fig:total_gen_vs_battery_cap}
\end{figure}

Figure~\ref{fig:solution} (right) shows that the total operational cost and generation costs first slightly decrease then increase when we increase the battery capacity on the transit buses. This suggests that the presence of the batteries can alleviate some strain on the system even though the total generation amount increases. Moreover, we can observe that the charging/discharging cost is mostly negative and decreasing up to a certain level. Even though operational constraints on the transit fleet require more energy to be stored, the price difference over the course of the day causes the transit authority to be able to arbitrage and gain some profit. 
In Figure \ref{fig:total_gen_vs_battery_cap}, we show that the generation profile changes dramatically as the battery capacities increase. Depending on the transit schedule, large batteries can cause added strain on the power system, while acting as generators during other periods.

\subsection{Case study for the stochastic formulations}
\label{sec:results_stoch}
\begin{figure}
\begin{center}
\includegraphics[trim= 0cm 0cm 0cm -0.3cm, width=8cm]{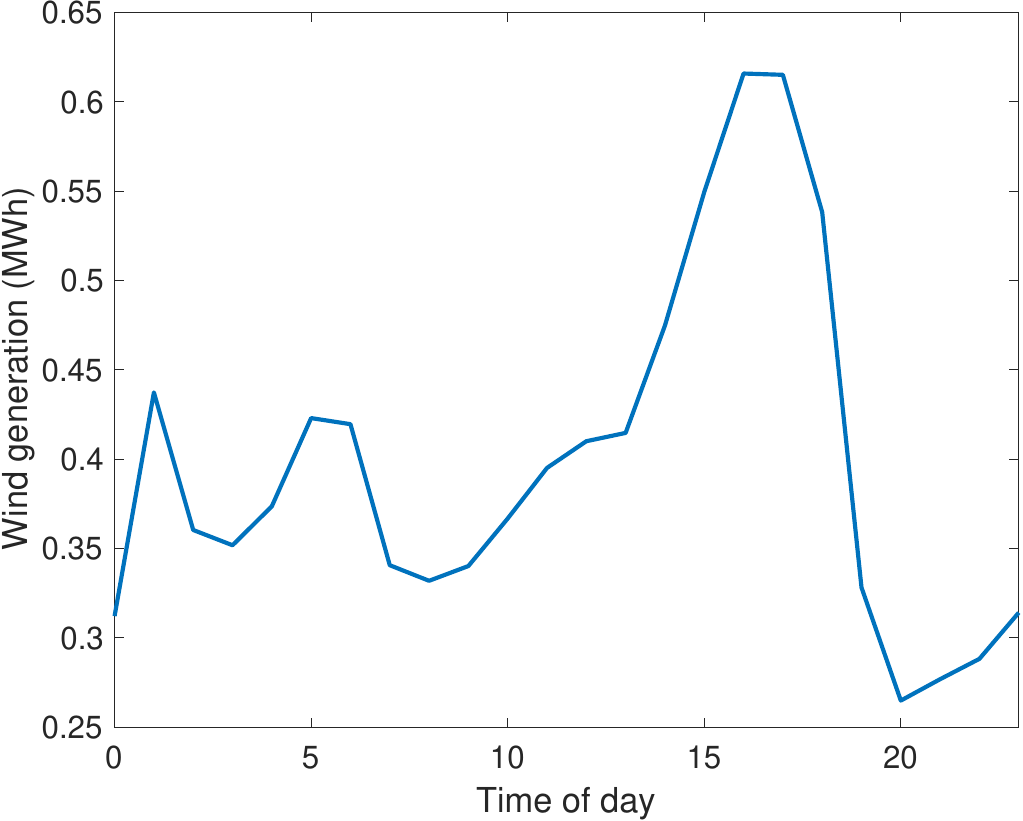}
\end{center}
\caption{Daily wind generation amount for 09/09/19 in San Jose (simulated)}
\label{fig:daily-wind}
\end{figure}

\begin{figure}
\begin{center}
\includegraphics[width=10cm]{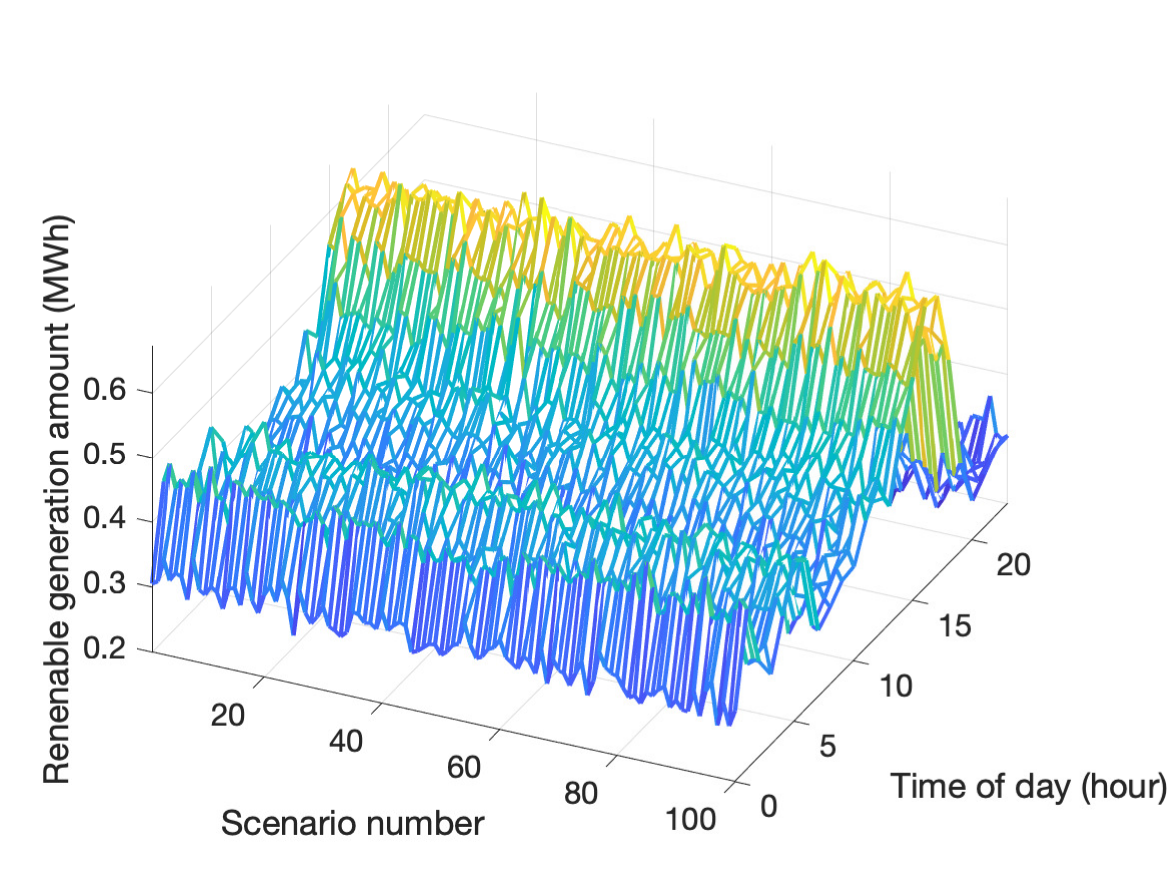}
\end{center}
\caption{Wind generation distribution over scenarios}
\label{fig:wind-sen-data}
\end{figure}

In this section, we consider the same instance described in Section~\ref{sec:results_det},  modifying the 9-bus system to include a wind generation unit located at node 4 in the power network displayed in Figure~\ref{fig:locations}. No other changes are made to the 9-bus system. The generation cost of wind is assumed to be the minimum over the linear cost coefficients of conventional generation, that is,
$$c^r_{it} = \min_{i \in \Ncal_g} \{ c^g_{it} \}.$$
In this manner, we ensure that the wind generator is always the cheapest among all generators, conventional or renewable. The wind generation data is obtained from a simulation API (see \cite{PFENNINGER20161251,STAFFELL20161224} for more details) for the day 09/09/19 within close proximity of the San Jose region with a maximum generation capacity set at 1 MWh. That is, the data are synthetic, but closely approximate the true wind generation of the area under consideration. Figure~\ref{fig:daily-wind} illustrates the daily data. Moreover, we add normally distributed noise with mean zero and a small variance (on the order of approximately 0.01 MWh) to create the desired number of scenarios. When 100 scenarios are considered, we obtain the distribution shown in Figure~\ref{fig:wind-sen-data}. 

The generated wind data are used as the renewable generation amounts for each scenario in both of the 2SSP formulations. The charging/discharging prices in the 2SSP formulations are derived as described in Section~\ref{sec:benefit_stoc}, and a visualization is provided in Figure~\ref{fig:prices_stoc}. In each of the following sections, we are limited to only 10 scenarios for renewable generation due to computational limitations. Note that the primary goal of these two models is to assess the capabilities of two different formulations in terms of wind utilization, since in reality, the renewable generation cost should be the smallest among all other costs, including the ramping cost of conventional generators. Hence, one should try to utilize wind generation as much as possible.

There is evidence of asymmetry in the ramping up and down of thermal generators. There are three identified sources for these costs: (1) {\em creep}, when components operate above the design temperature; (2) {\em thermal fatigue}, when changes in temperature result in mechanical failure; (3) {\em creep--fatigue} interactions, when the two effects above compound. As studied by \cite{MOAREFDOOST-ramp-cost}, in the ramp-up of thermal generators, these three effects are present, whereas in ramp-down the main effect is thermal fatigue. For this reason, in the following numerical results, we use ramping costs defined by:
\begin{align*}
    c_{it}^{g,+} &= 1.2  c_{it}^{g} \\
    c_{it}^{g,-} &= 0.5  c_{it}^{g}.
\end{align*}
That is, for ramping costs we only consider the linear generation cost coefficient $c_{it}^{g}$ to make ramping slightly cheaper than generation, and, ramping down is much cheaper than ramping up.

\subsubsection{Case study for the ramping-based formulation}

\begin{figure}
\begin{center}
\includegraphics[trim = 0cm 0cm 0cm -0.5cm, width=10cm]{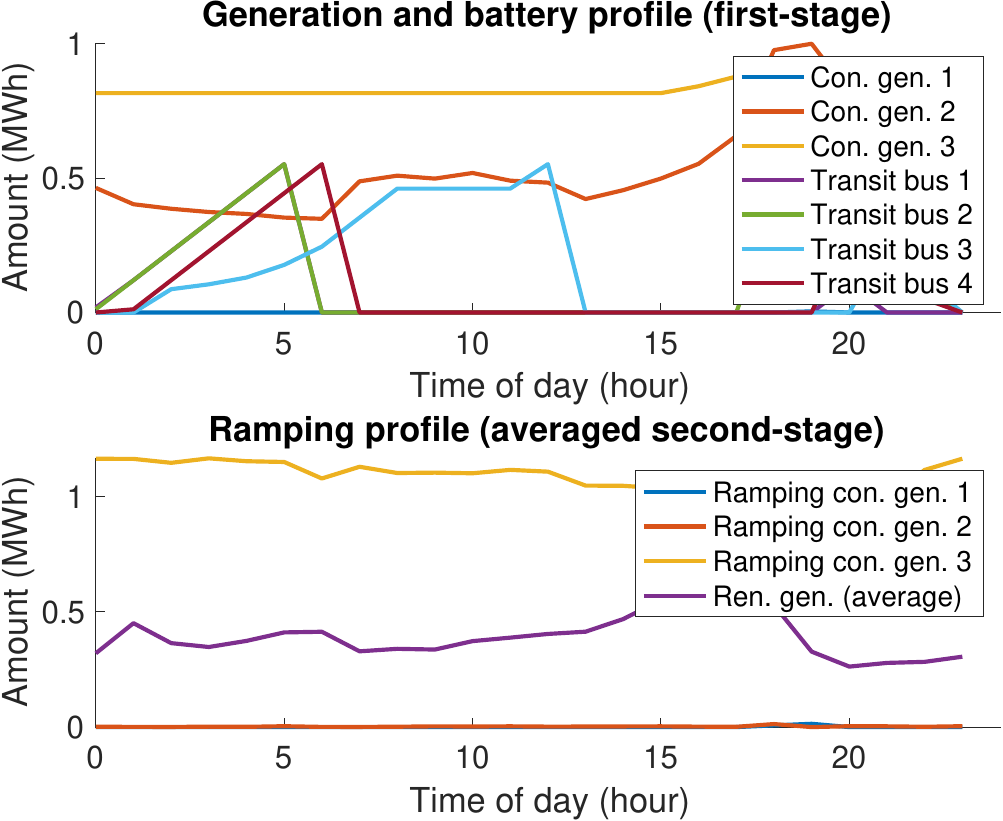}
\end{center}
\caption{Optimum generation, charging, and ramping levels }
\label{fig:solution_stoch_ramping}
\end{figure}

Figure~\ref{fig:solution_stoch_ramping} gives a summary of the optimal solution of the formulation presented in Section~\ref{sec:stoc_model_ramp}. We observe that only one generator uses ramp-up in the second stage.
Note that the behavior of the battery levels is similar to that resulting from the solution to the deterministic formulation, in that they have a similar charging pattern; they charge in the early hours of the morning, as well as late at night. That is, regardless of whether or not the decision is being made by the social planner, the risk associated with deciding when to charge/discharge does not drastically affect the charging/discharging times of the transit buses. This can be attributed to the fact that the generators are flexible enough to handle fluctuations in the amount of renewable generation.

\begin{figure}
\begin{center}
\includegraphics[width=10cm]{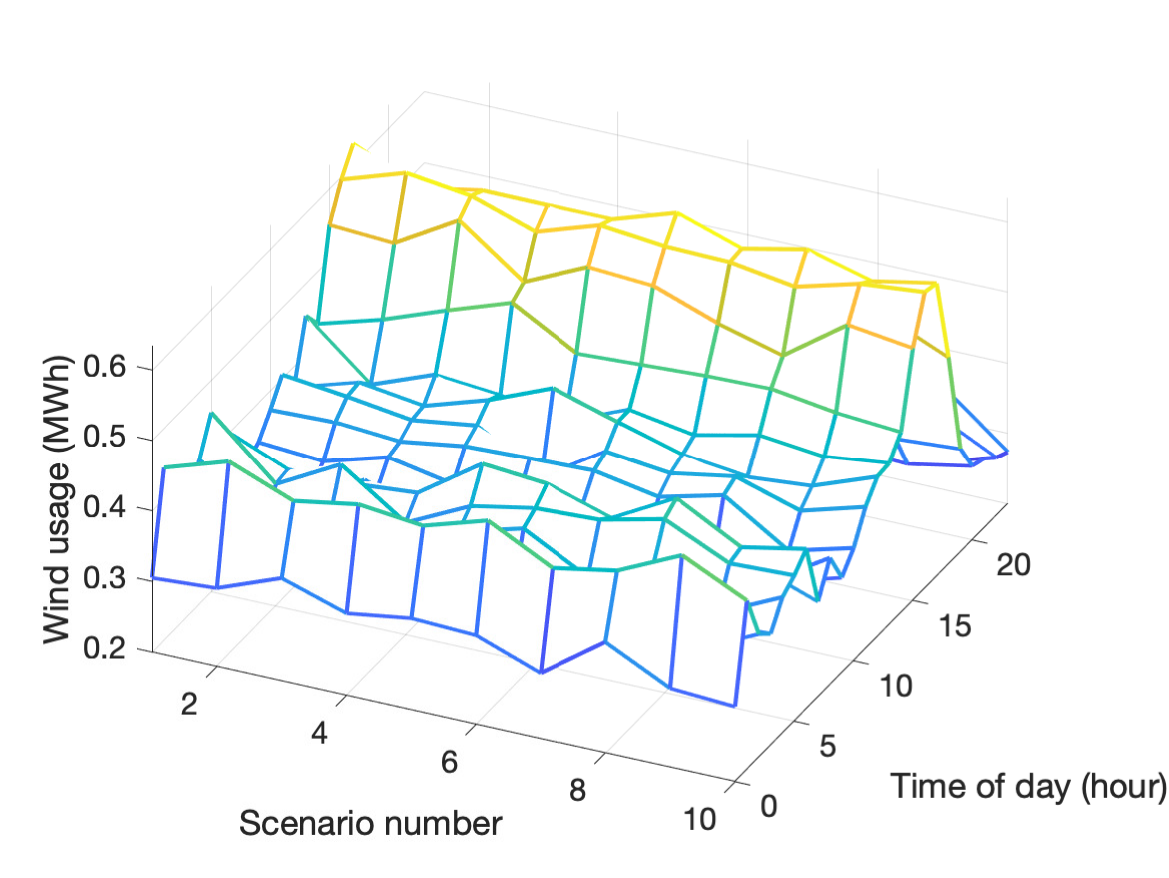}
\end{center}
\caption{Wind usage amount in the optimal solution}
\label{fig:wind_usage_stoch_ramping}
\end{figure}

In Figure \ref{fig:wind_usage_stoch_ramping}, we present the optimal values of wind generation corresponding to the second stage resulting from solving the optimization problem found in Section \ref{sec:stoc_model_ramp} across all time periods and scenarios. For each scenario, we find that wind-generated energy is utilized in its entirety. We calculate the wind utilization measure as follows:
\begin{equation}\label{e:windUt}
    \text{Wind Utilization} = \frac{\text{total wind usage in optimal solution}}{\text{aggregation of wind generation from scenarios}}.
\end{equation}
In essence, wind utilization calculates the ratio of utilized wind generation (in the optimal solution) to the total available wind generation over all of the scenarios. Since we consider all of the scenarios in the calculation of the measure, it can be thought of as an average utilization over the scenarios. This is useful to quantify in general, since we would like to utilize as much wind energy as possible in the solution.

In our numerical results corresponding to the optimization problem found in Section \ref{sec:stoc_model_ramp}, the wind utilization is found to be 1.0 (i.e., 100\%). As such, in this case, ramping as a recourse action is flexible enough to handle the deviations in wind generation over the entire planning horizon. However, if we restrict the ramping quantity in the second stage to only allow for very small deviations in the generation, this total utilization might not be attainable. This is because the wind generator is limited to a total hourly (periodic) generation capacity of 1 MWh, whereas in our numerical experiments we allow for ramping that can introduce fluctuations on the order of approximately 100 MWh. 

\subsubsection{Case study for the charging/discharging-based formulation}

\begin{figure}
\begin{center}
\includegraphics[width=10cm]{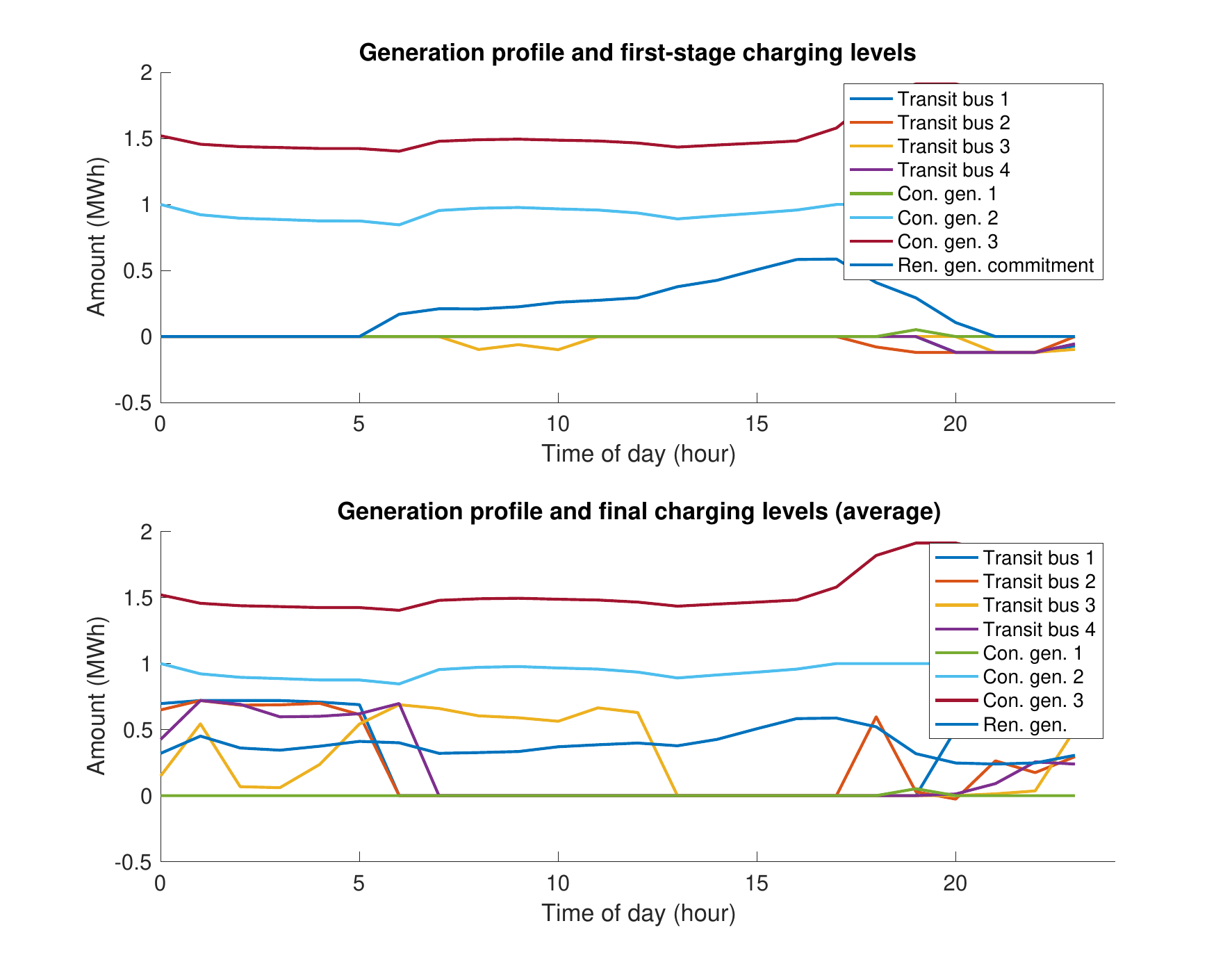}
\end{center}
\caption{Optimum generation and charging levels}
\label{fig:solution_stoch_charging}
\end{figure}

In Figure~\ref{fig:solution_stoch_charging}, we provide a summary of the solution of the formulation described in Section~\ref{sec:stoc_model_charge}. It is clear that the generator outputs are much higher in the first stage compared to the solution of the alternative model displayed in Figure~\ref{fig:solution_stoch_ramping}, as second-stage adjustment via ramping is not allowed in this formulation. Moreover, we also observe that the transit charging/discharging levels change in the second stage rather than in the first stage, to utilize the wind energy as much as possible. 

\begin{figure}
\begin{center}
\includegraphics[width=15cm]{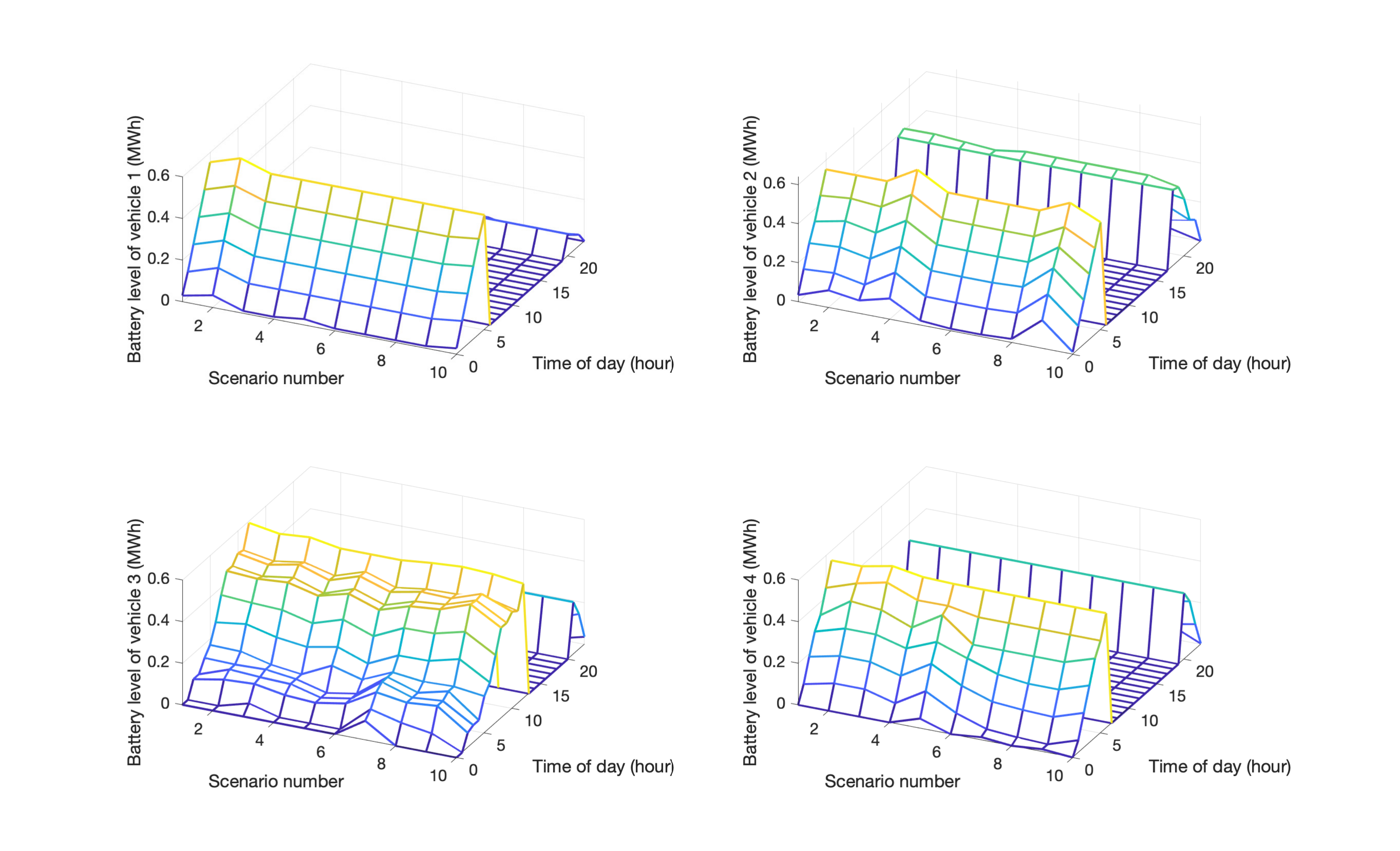}
\end{center}
\caption{Battery levels of vehicles}
\label{fig:soc_stoch_charging}
\end{figure}

Additionally, in Figure~\ref{fig:soc_stoch_charging} one can view the battery levels of transit vehicles for each scenario (for 10 scenarios, as mentioned above). We can observe that the entire fleet adjusts its battery levels based on the availability of renewable generation while ensuring that the operational constraints are feasible. In particular, these operational constraints limit the solution space and make the recourse actions much more difficult than those of simple batteries.

\begin{figure}
\begin{center}
\includegraphics[width=10cm]{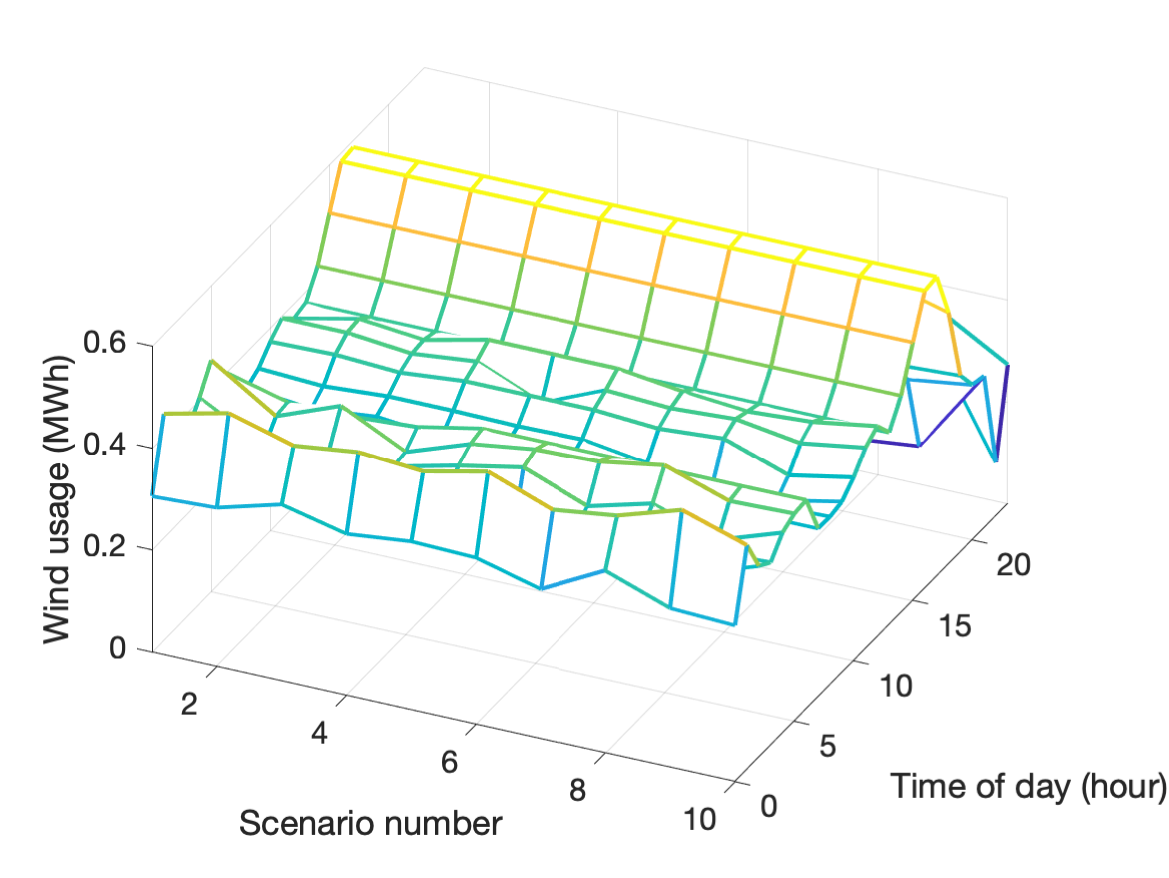}
\end{center}
\caption{Wind usage amount in the optimal solution}
\label{fig:wind_usage_stoch_charge}
\end{figure}

Similar to the previous subsection, we calculate the wind utilization measure given by equation \eqref{e:windUt} and find that for the solution to the optimization problem given in Section \ref{sec:stoc_model_charge} the total wind utilization is 0.9645 (96.45\%). In Figure \ref{fig:wind_usage_stoch_charge}, we present the optimal values of wind generation corresponding to the second stage resulting from solving the the optimization problem found in Section \ref{sec:stoc_model_charge} across all time periods and scenarios. One observation is that in this case, the solution is limited by operational constraints and battery capacities. As a result, compared to the previous utilization displayed in Figure~\ref{fig:wind_usage_stoch_ramping}, we find that the model from Section \ref{sec:stoc_model_charge} does not fully utilize wind generation later in the day (roughly from $t=15$ on). These underutilized time periods occur specifically when the fleet is serving public transit demand. However, if the fleet contains numerous transit buses with complementary off-schedules that span an entire day, the utilization would certainly be superior.

\subsection{Larger power grids}

\label{diff-power-grid-results-det}
To verify the capabilities of the proposed formulation for cooperation, we further test it when larger power networks are coupled with the transit network. In particular, we utilize the \textsc{MATPOWER} case files \citep{MATPOWER}, specifically, \texttt{case14}, \texttt{case30}, \texttt{case39}, \texttt{case57}, and \texttt{case118}, where the number denotes the number of nodes in the system. These cases are directly accessible within the \textsc{MATPOWER} instance list. We only modify these instances by scaling down the demand, which is done similarly to the procedure explained in Section~\ref{sec:results_det}. Note that the information associated with the transit fleet remains the same. We always assume the charging stations are connected at the first six nodes of the power system.

We obtain similar findings in these instances. While imposing additional demand on the power network, the transit fleet can act as batteries for some periods and alleviate generation costs by smoothing the generation profile. Moreover, in instances in which congestion is present, we observe the effect of the transit fleet more clearly, since the unit prices of electricity set by the optimal dual variables associated with the nodal balance equations can dramatically change over time (due to congestion), and the transit fleet can utilize an arbitrage strategy, which, in fact, benefits the whole system in cooperation. 

Table~\ref{tab:sol-times-det} provides a summary of solution times for the deterministic formulation. The solution time of the deterministic formulation with our modified \texttt{case9} instance is much higher than the rest of the cases. This is primarily related to the additional congestion introduced into the power network. Aside from the \texttt{case9} outlier, we observe that the solution times increase when we incorporate larger power networks.

\begin{table}
 \caption{Solution times of the deterministic formulation}
 \label{tab:sol-times-det}
 \begin{center}
 \begin{tabular}{lr}
    \toprule
      \textbf{Instance} & \textbf{Solution time (s)} \\
      \midrule
      \texttt{case9} & 23.1256 \\
      \texttt{case14}  	& 0.3720\\
      \texttt{case30}  	& 0.4434\\
      \texttt{case39} 	& 0.4810\\
      \texttt{case57}  	& 0.4285\\
      \texttt{case118}  & 1.1776\\
      \bottomrule
    \end{tabular}
\end{center}
\end{table}

Table~\ref{tab:summary-stoch-ramp} provides a summary of results for the ramping-based stochastic formulation on larger power networks. From the table, one can observe that we have complete wind utilization in every case because ramping limits of conventional generators are large enough to compensate for the uncertainty in the second stage. Moreover, in general, solution times increase when we integrate a larger power network into the formulation.

\begin{table}
\caption{Summary of results for the ramping-based formulation}
\label{tab:summary-stoch-ramp}
\begin{center}
\begin{tabular}{lrr}
    \toprule
      \textbf{Instance} & \textbf{Solution time (s)} & \textbf{Wind utilization} \\
    \midrule
      \texttt{case9}	& 0.8120 & 1\\
      \texttt{case14}	 & 2.3466 & 1\\
      \texttt{case30} & 4.9761 & 1\\
      \texttt{case39} 	 & 6.8127 & 1\\
      \texttt{case57}	 & 9.4370 & 1 \\
      \texttt{case118} & 146.8926 & 1 \\
     \bottomrule
    \end{tabular}
\end{center}
\end{table}

Finally, Table~\ref{tab:summary-stoch-charge} provides a summary of results for the charging/discharging-based stochastic formulation on larger power networks. We can observe that in all of the cases, we have the same wind utilization (0.9645), since the charging limit of the transit fleet remains the same in all of the formulations, and the fleet can compensate for a portion of the uncertainty in the second stage. More importantly, we also observe that allowing second-stage charging/discharging of the fleet introduces additional complexity to the formulation since solution times are much larger than those in Table~\ref{tab:summary-stoch-ramp}.

\begin{table}
\caption{Summary of results for the charging/discharging-based formulation}
\label{tab:summary-stoch-charge}
\begin{center}
    \begin{tabular}{lrr}
    \toprule
      \textbf{Instance} & \textbf{Solution time (s)} & \textbf{Wind utilization} \\
      \midrule
      \texttt{case9} & 121.3390 & 0.9645\\
      \texttt{case14} & 1690.1862 & 0.9645\\
      \texttt{case30} & 32.7924 & 0.9645\\
      \texttt{case39} & 47.3546 & 0.9645\\
      \texttt{case57} & 52.4875 & 0.9645 \\
      \texttt{case118} & 991.3013 & 0.9645 \\
      \bottomrule
    \end{tabular}
\end{center}
\end{table}

\section{Managerial insights}
\label{sec:managerial_insights}

\subsection{Benefit of coordinated optimization} \label{sec:benefit_det}
In this section, we use the deterministic formulation to evaluate the benefit of employing a coordinated strategy between the ISO and the transportation authority. To serve as a baseline, we consider their uncoordinated operation where each party acts to manage its own objective. The uncoordinated optimization scheme is summarized in Procedure~\ref{uncoord_opt_scheme}.

\begin{algorithm}
\caption{Uncoordinated optimization scheme}
\label{uncoord_opt_scheme}
\begin{algorithmic}[1]
\State Feasible charging scenarios are generated by ISO \label{alg:sen} 
\For{each scenario}
\State ISO solves dispatch problem with additional demand determined by the scenario and obtains a set of prices from dual variables associated with nodal balance constraints \label{alg:mpopf} 
\State Using the prices, transit authority optimizes its own problem and obtains a charging/discharging policy \label{alg:charge} 
\State  Charging/discharging policy is realized by ISO, and ISO objective value is obtained \label{alg:eval_ISO}
\State  Charging/discharging policy is evaluated under baseline prices to obtain transit objective
\State Total uncoordinated cost is calculated \label{alg:eval_transit}
\EndFor
\end{algorithmic}
\end{algorithm}

Step~\ref{alg:sen} of Procedure~\ref{uncoord_opt_scheme} generates a given number of scenarios, which only takes feasible charging of batteries into consideration. That is, the charging anticipated by the ISO is guaranteed to satisfy transit-operational constraints. This proves to be beneficial to the ISO, as this procedure can be seen as educated anticipation from the ISO, where they have access to some information related to the transportation system. Step~\ref{alg:mpopf} simply solves an MPOPF without any of the transportation aspects in the formulation described in Section~\ref{sec:model_det}. In Step~\ref{alg:charge}, only the transportation aspects such as charging/discharging and location/relocation of the transit buses are considered in a separate formulation. Then, Step~\ref{alg:eval_ISO} solves the formulation in Step~\ref{alg:mpopf} with the optimal charging/discharging obtained in Step~\ref{alg:charge} as an additional demand. Next, to make the comparison fair, in Step~\ref{alg:eval_transit}, we consider a set of baseline prices to evaluate the solution obtained in Step~\ref{alg:charge}. This can also be seen as the average prices over the scenarios anticipated by the ISO.

For coordinated optimization, the deterministic formulation in Section~\ref{sec:model_det} is solved, with baseline prices and the objective value being calculated with convex combination coefficient $\alpha = 0.5$ in the objective function \eqref{eq:obj1}. Then, since the uncoordinated cost accounts for the summation of the two costs (rather than a convex combination), we scale down the total uncoordinated cost by halving it, to make a fair comparison.

\begin{figure}
\begin{center}
\includegraphics[trim= 0cm 0cm 0cm -0.25cm, width=10cm]{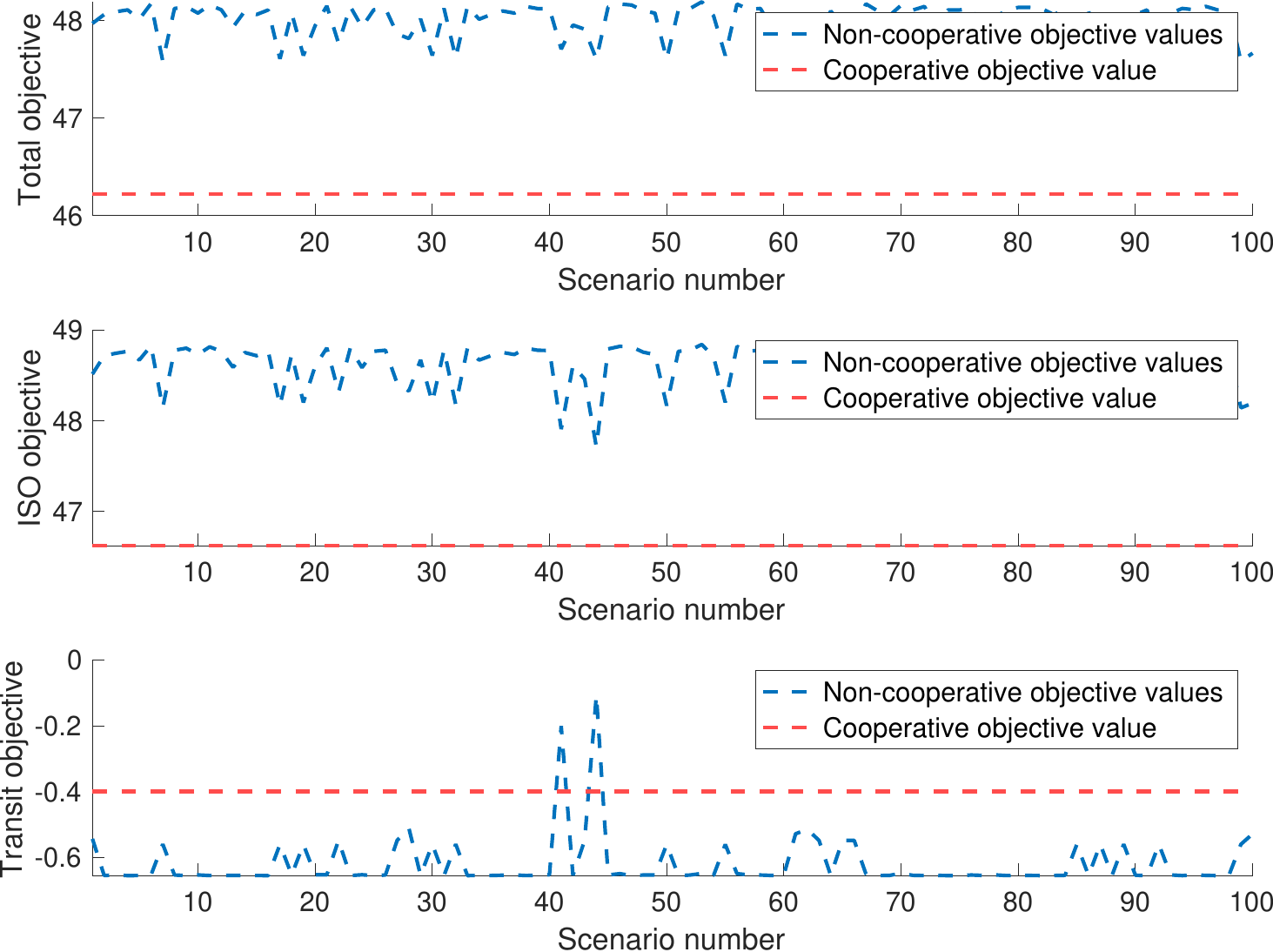}
\end{center}
\caption{Comparison of objective values}
\label{fig:benefit_analysis}
\end{figure}

Figure~\ref{fig:benefit_analysis} details the uncoordinated objective values for the transit authority, ISO, and their sum (all of the quantities are scaled down by $0.5$ due to the presence of $\alpha$ in the co-optimization), and compares these with the coordinated objective values. It is immediate to see that the total and ISO objective values are always worse in the non-cooperative case. An immediate reason is that the anticipation of charging schedules made by the ISO is quite different from the actual strategy employed by the transportation authority, and the power system incurs additional generation cost. However, since we consider charging prices that are averaged over scenarios as the baseline prices, the transit objective values are comparable between the two strategies. These results suggest that the ISO benefits largely from cooperation with the transit authority. This was expected since the ISO minimizes its cost when it has full knowledge of the optimal charging schedule, and anticipation of any deviation from the optimal charging schedule will of course worsen the ISO's operation. 

We also tested the benefit of coordination in the larger power networks discussed in  Section~\ref{diff-power-grid-results-det}. Table~\ref{tab:benefit-summary-det} summarizes the results. Note that the non-cooperative objectives are calculated by taking the average of the objectives over 100 different anticipation scenarios. As expected, on average, the transit objectives in both strategies match, since we use the average prices (over 100 anticipation scenarios) as our electricity prices in the evaluation step. Moreover, we observe that in all of the cases the ISO benefits from the cooperative strategy, with the magnitude of the benefit becoming less significant when larger power networks are considered. This supports our previous results concluding that the scale of the battery capacities (or the transit fleet) plays a crucial role in the benefit obtained by a cooperative strategy.
Observing the intricacies of the benefit analysis on the stochastic formulations is yet another exploration. However, for the sake of conciseness, we believe that the primary aspect of coordination is demonstrated.

\begin{table}
\caption{Benefit analysis of the deterministic formulation (summary)}
\label{tab:benefit-summary-det}
\begin{center}
\begin{tabular}{lrrrr}
    \toprule
      \textbf{Instance} & \textbf{Coop. t.$^a$ obj.} & \textbf{Non-coop. t.$^a$ obj.} & \textbf{Coop. I.$^b$  obj.}  & \textbf{Non-coop. I.$^b$ obj. } \\
      \midrule
      \texttt{case14} & 17.7027 & 17.7028 & 752.0961 & 752.0983\\
      \texttt{case30} & 1.5947 & 1.5950 & 126.2755 & 127.8705\\
      \texttt{case39} & 0.2724 & 0.2724 & 23.5928 & 23.5929\\
      \texttt{case57} & 17.5882 & 17.5882 & 900.5379 & 900.5385\\
      \texttt{case118} & 17.5838 & 17.5838 & 4092.9120 & 4092.9150\\
      \texttt{case145} & 17.5269 & 17.5269 & 1320.0214 & 1320.0215\\
      \bottomrule
\end{tabular}
\end{center}
\footnotesize{$^a$ Transit authority, $^b$ ISO}
\end{table}

\subsection{Potential benefits of alternative recourse actions} \label{sec:benefit_stoc}
In this section, we are interested in a comparison of the two stochastic formulations presented in Sections~\ref{sec:stoc_model_ramp} and \ref{sec:stoc_model_charge}. Recall that the first model investigates recourse actions taken by the ISO (ramping), whereas the second model considers the actions taken by the transit operator (using batteries to handle the randomness in the second stage). 

\begin{figure}
\begin{center}
\includegraphics[width=10cm]{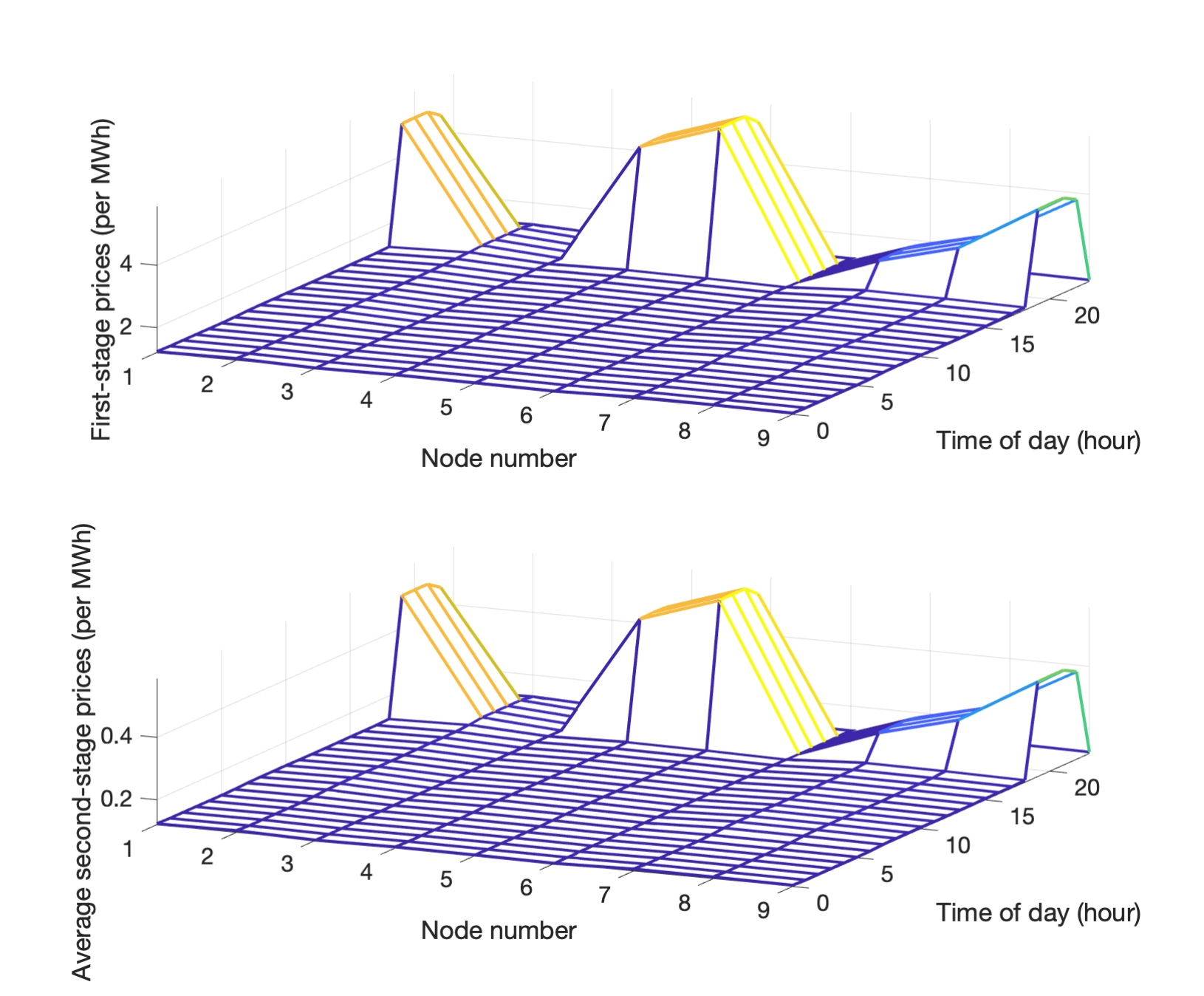}
\end{center}
\caption{First-stage and second-stage prices for charging/discharging of the transit fleet}
\label{fig:prices_stoc}
\end{figure}

To ensure a fair comparison, we consider the same set of charging/discharging prices in both of the formulations. Note that in the second formulation, there are two sets of prices, corresponding to the prices in the first stage, $c_{it}$, and in the second stage, $c_{it\omega}^+$. These prices are obtained by solving the two-stage stochastic multi-period dispatch problem given in Appendix~\ref{apx:MPOPF}. In more detail, the optimal values for the dual variables associated with the first-stage nodal balance equations determine the first-stage prices, and the second-stage prices are derived similarly. Figure~\ref{fig:prices_stoc} illustrates the first-stage prices and the average second-stage prices for charging/discharging of the transit fleet.

In both stages, we observe that prices sharply increase near the end of the day due to the extra stress imposed by the extra wind generation around similar times as previously shown in Figure~\ref{fig:daily-wind}. Moreover, the second-stage prices are much lower than the first-stage prices. This could be attributed to the fact that the second-stage flows are much smaller, as the scale of recourse actions is typically much smaller than that of the first-stage decisions.

Next, we analyze the costs in each of the two models provided in Sections~\ref{sec:stoc_model_ramp} and \ref{sec:stoc_model_charge}. The following  parameters are relevant:
\begin{itemize}
    \item First-stage and second-stage charging/discharging prices $(c_{it}, c_{it\omega}^+)$
    \item Renewable generation cost $(c_{it}^r)$
    \item Ramp-up/down costs for conventional generators $(c_{it}^{g,+}, c_{it}^{g,-})$
\end{itemize}

Observe that the renewable generation term appears in both of the objectives \eqref{eq:obj1_s1} and \eqref{eq:obj1_s2}. Then, the only freedom we have is in determining ramping costs for conventional generators. Based on this, we vary the ramping cost values and obtain a trade-off between the costs of the two models. We provide the scheme for comparing the two stochastic formulations in Procedure~\ref{alg:benefit_analysis_stochastic}.

\begin{algorithm}
\caption{Efficacy analysis for the two stochastic formulations}
\linespread{0.6}\selectfont
\label{alg:benefit_analysis_stochastic}
\begin{algorithmic}[1]
\For{each ramp-up/down cost}
\State \label{alg:mpopf_stoch} Using ramping costs, solve the two-stage stochastic MPOPF to obtain first-stage and second-stage charging/discharging prices
\State \label{alg:ramping_stoch} Using the first-stage prices and ramping costs, solve the model with ramping as the recourse action and obtain an objective value
\State \label{alg:charging_stoch} Using the first-stage and second-stage prices, solve the model with charging/discharging as the recourse action and obtain an objective value
\EndFor
\end{algorithmic}
\end{algorithm}

\begin{figure}
\begin{center}
\includegraphics[width = 10cm]{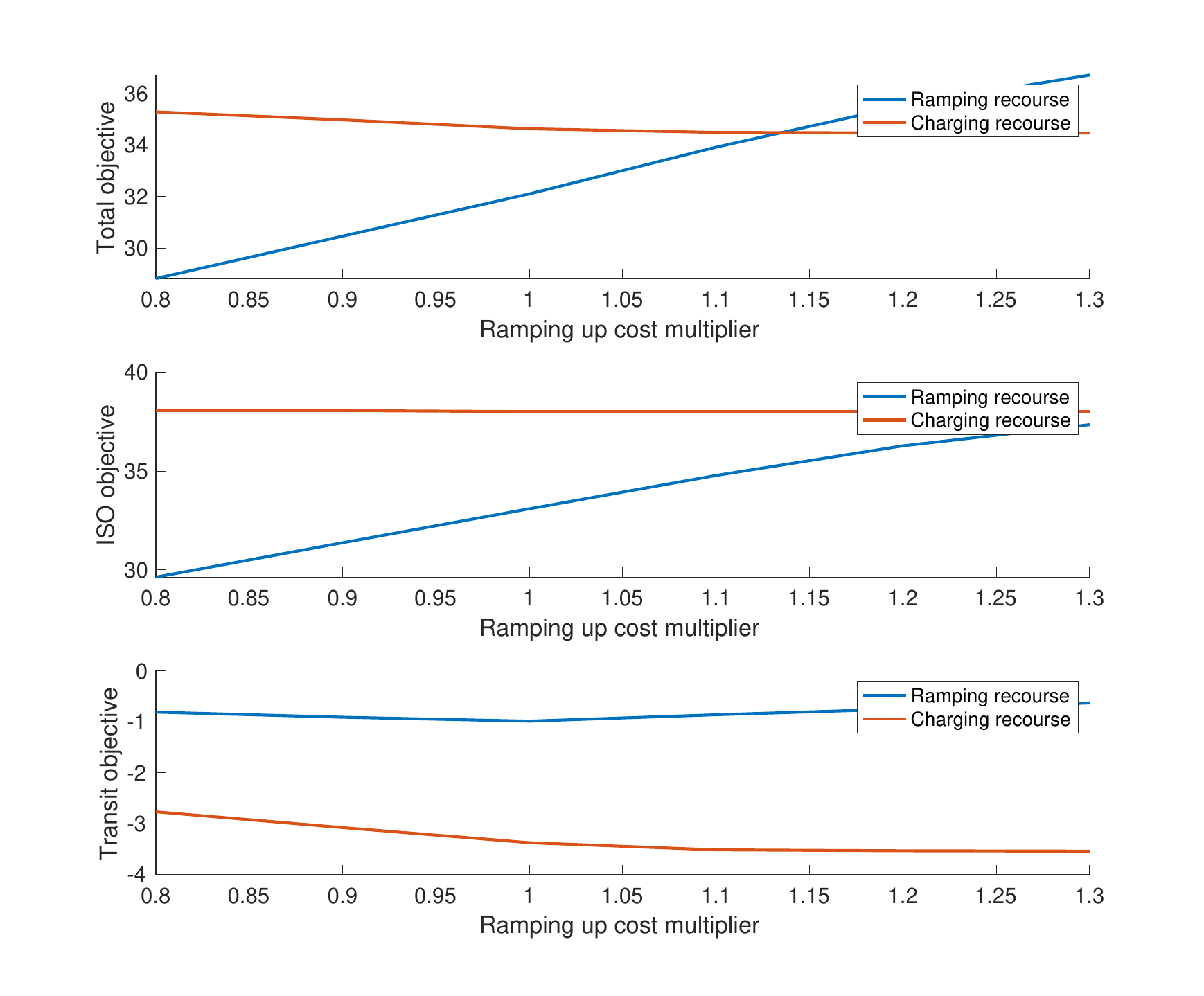}
\end{center}
\caption{Comparison of objective values in stochastic formulations by varying ramp-up cost}
\label{fig:comp_obj_stoch}
\end{figure}

In Procedure~\ref{alg:benefit_analysis_stochastic}, we provide the scheme to obtain the trade-off by calibrating the ramping costs in objective \eqref{eq:obj1_s1}. In more detail, we systematically vary the ramping costs and obtain the adjusted prices for both first and second stages via locational marginal prices in step~\ref{alg:mpopf_stoch}. Currently, for simplicity, in step~\ref{alg:mpopf_stoch}, we do not consider the demand added by the transit fleet. Alternatively, one could also incorporate an average demand for transit charging/discharging. Then, by way of these new prices, we separately solve the two formulations and compare their objective values. Note that, in our experiments, we only vary the ramp-up cost by changing the multiplier $\gamma$ in the following manner:
$$
    c_{it}^{g,+} = \gamma  c_{it}^{g} ~\text{and}~
    c_{it}^{g,-} = 0.5  c_{it}^{g},
$$
where $\gamma \in \{0.8,0.9,1,1.1,1.2,1.3 \}$. By allowing the value of $\gamma$ to vary over these values, we obtain the results in Figure~\ref{fig:comp_obj_stoch}.

As can be seen, when the ramp-up cost multiplier value is around 1, the total objective values in the two models become competitive. Since having the value around 1 is reasonable, we can conclude that the recourse action provided by charging, as opposed to ramping, can be a useful alternative depending on the generator ramping costs. Moreover, it is interesting to observe that the transit objective in the second formulation is much smaller. One possible reason for this could be that the batteries have more flexibility in the second formulation since they can also arbitrage between the first and second stages, whereas in the first formulation, they can only charge/discharge in the first stage.

We further conducted a similar analysis when larger power networks are coupled with the transit network. In general, we obtained similar findings in which the trade-off values for $\gamma$ were close to 1. This supports the idea that the two 2SSP models can be alternatives to each other depending on ramping costs, and more importantly, total battery capacity. We omit further details to avoid displaying repetitive findings.

In a broad view, in this section, we compared two types of recourse actions in our framework. It is crucial to note that for additional flexibility promoting the integration of renewable generators, one can utilize these two recourse actions simultaneously.

\subsection{Pricing and co-optimization}
\label{sec:pricing-and-coopt}
In this section, motivated by the work of \cite{kok-pricing-impact-on-inv}, which analyzes the effect of pricing policies (flat pricing vs.\ peak pricing) on the investment levels of renewable and conventional sources from the perspective of utility firms, we investigate co-optimization under different pricing policies for the transit authority. Specifically, we compare peak pricing (pricing obtained from locational marginal prices) to flat pricing (which is obtained by averaging peak prices) in different aspects.

\begin{figure}
\begin{center}
\begin{subfigure}{.43\textwidth}
  \includegraphics[width=7cm]{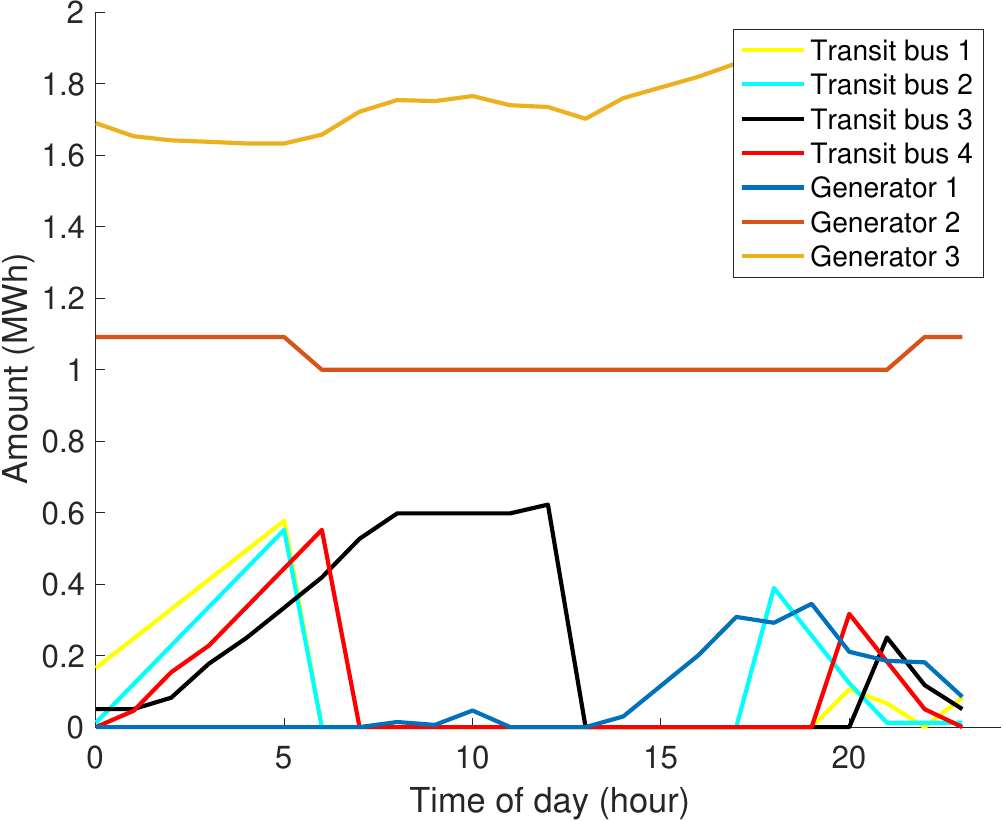}
\end{subfigure}
\begin{subfigure}{.43\textwidth}
    \includegraphics[width=7cm]{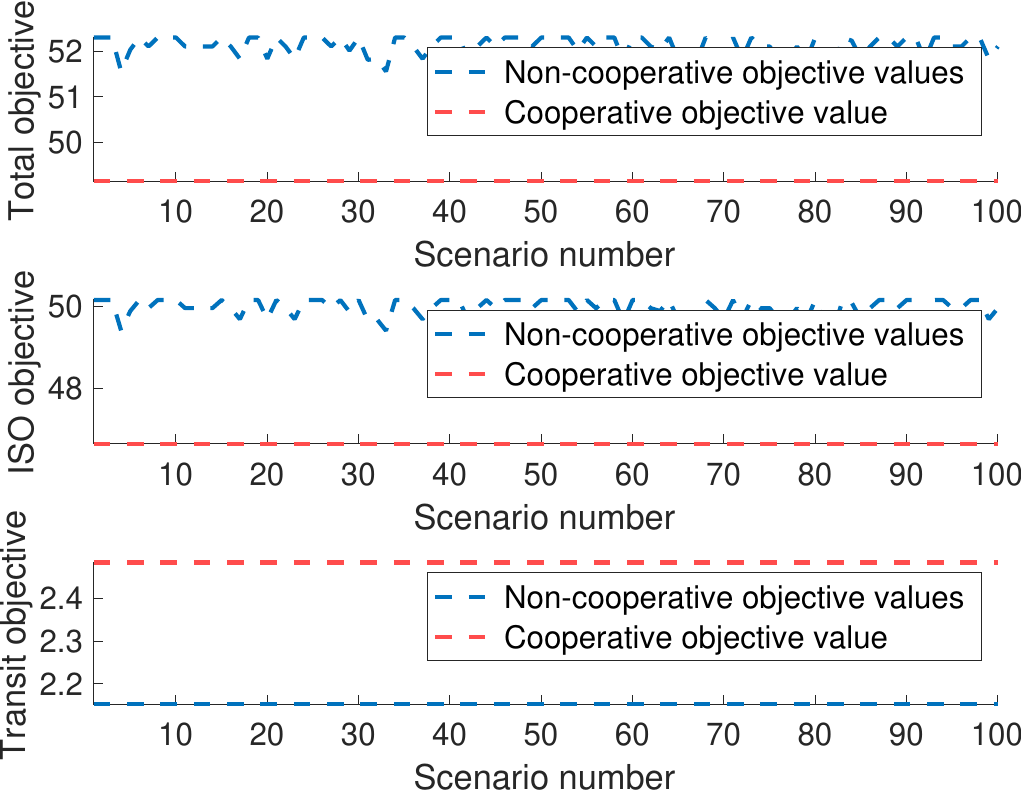}
\end{subfigure}
\end{center}
\caption{Deterministic co-optimization results under flat pricing: generation and battery level profile (left) and comparison of objective values (right)}
\label{fig:det_coopt_flat}
\end{figure}

In Figure~\ref{fig:det_coopt_flat}, we provide the optimal solution, along with the benefit analysis under flat pricing. If we compare Figure~\ref{fig:det_coopt_flat} (left) to the solution under peak pricing displayed in Figure~\ref{fig:solution} (left), we can observe that under flat pricing, both charging and generation fluctuate less, simply due to the elimination of arbitrage for the transit authority. However, since peak prices can be thought of as proxies for congestion information, the indirect benefit (in terms of relieving congestion) of utilizing peak prices is lost under flat pricing. Further, if we compare Figure~\ref{fig:det_coopt_flat} (right) to its counterpart in Figure~\ref{fig:benefit_analysis}, we immediately observe that costs increase under flat pricing. One immediate reason is that arbitrage is not possible under flat pricing and transit objectives increase dramatically. It is interesting to observe that under different scenarios, non-cooperative transit charging under flat pricing remains the same since buses only charge to satisfy their operational requirements. Moreover, cooperative transit objective values are higher than non-cooperative values because the fleet is charging more in total due to two reasons: firstly, the fleet is relocating to decrease the system cost; and secondly, the fleet is discharging to decrease generation cost and the efficiency of batteries is less than 1 ($\eta < 1$). Moreover, these compromises made by the transit authority save more for the ISO, thus the total system cost decreases significantly.

\begin{figure}
\begin{center}
 \begin{subfigure}{0.475\textwidth}
            \includegraphics[width=\textwidth]{figures_v2/wind_usage_charging_case9-eps-converted-to.pdf}
            \caption[]%
            {\footnotesize Under peak pricing in both of the stages} 
        \end{subfigure}
        \begin{subfigure}{0.475\textwidth}  
            \includegraphics[width=\textwidth]{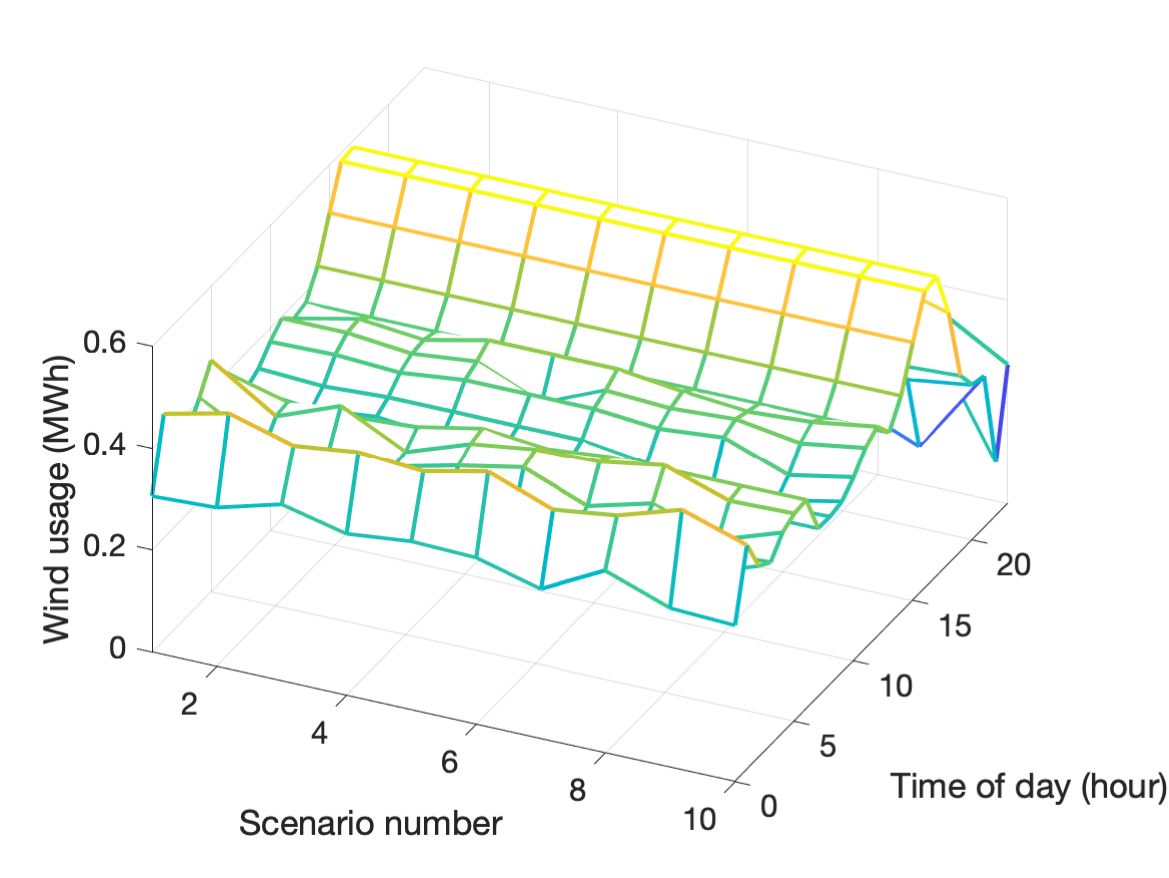}
            \caption[]%
            {\footnotesize Under flat pricing in the first stage}   
        \end{subfigure}
        \vskip -0.3cm
        \begin{subfigure}{0.475\textwidth}   
            \includegraphics[width=\textwidth]{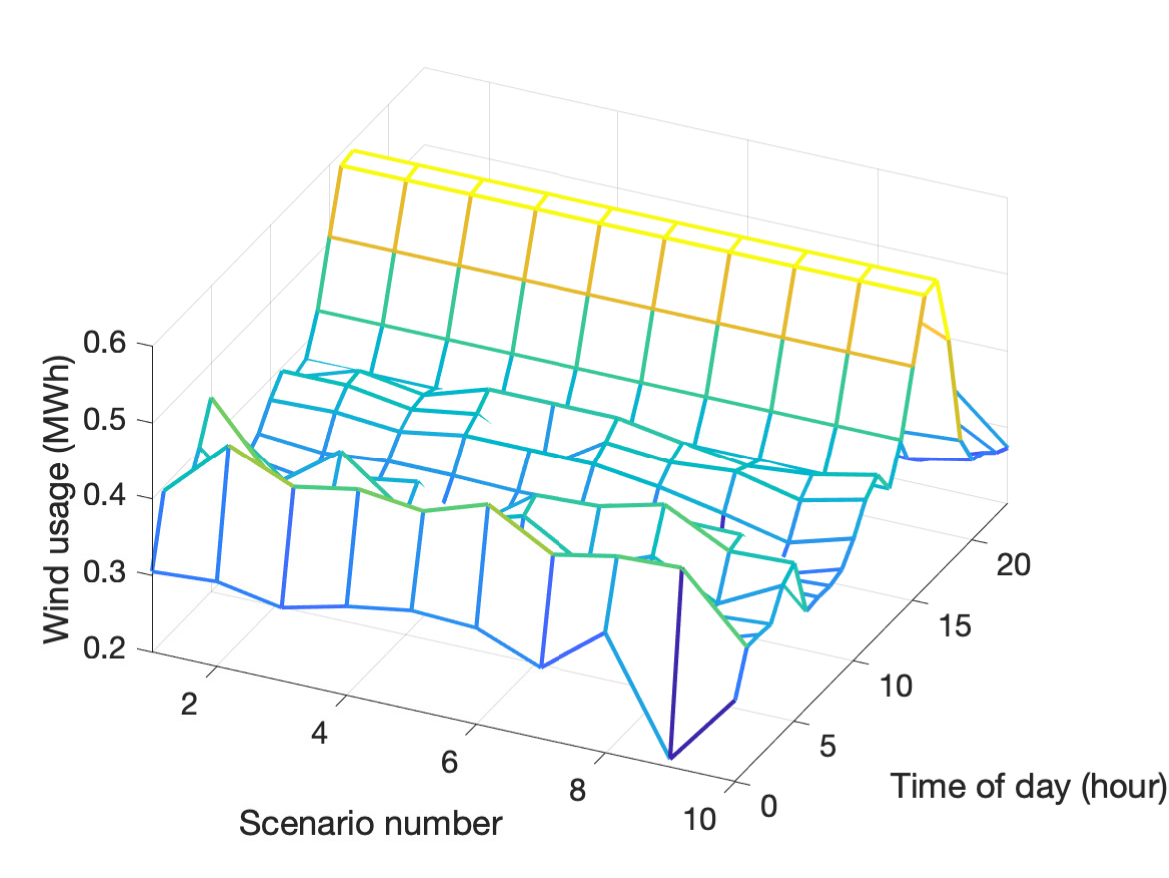}
            \caption[]%
            {\footnotesize Under flat pricing in the second stage}  
        \end{subfigure}
        \begin{subfigure}{0.475\textwidth}   
            \includegraphics[width=\textwidth]{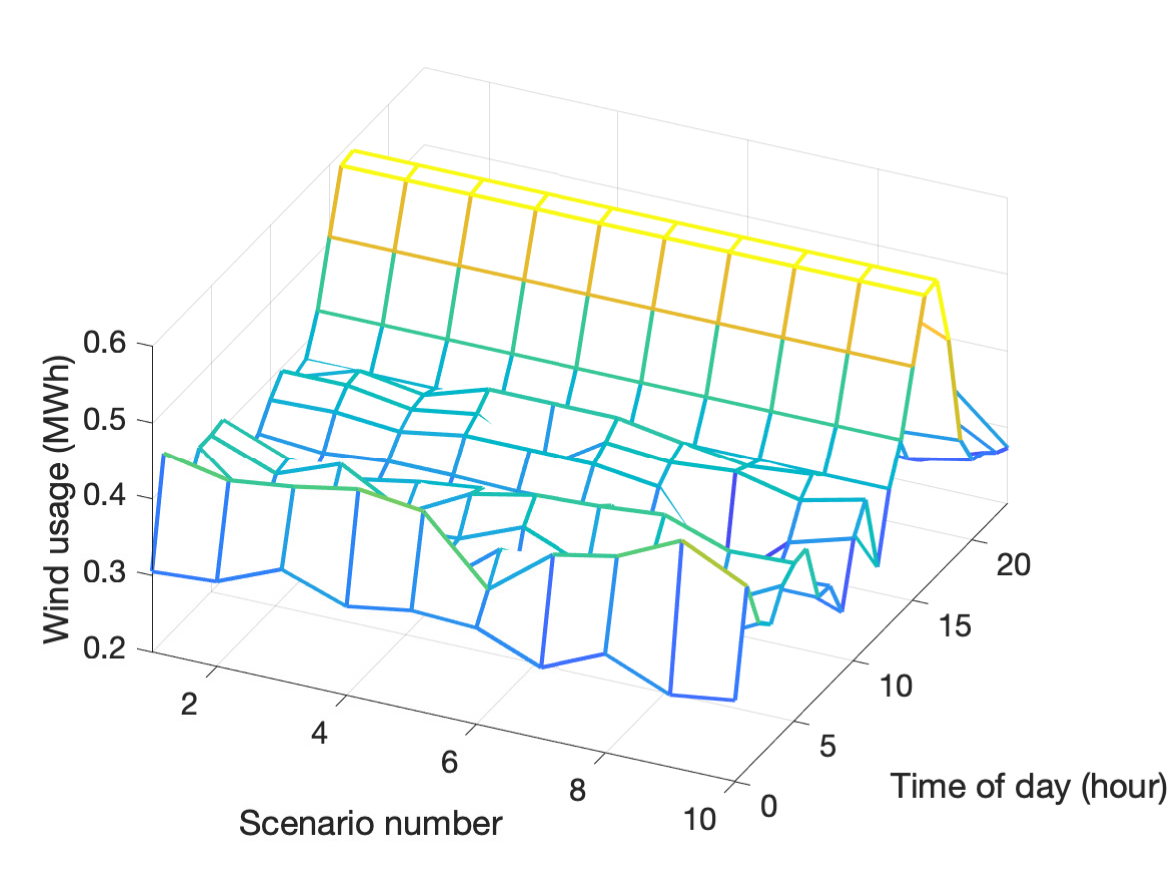}
            \caption[]%
            {\footnotesize Under flat pricing in both of the stages}
        \end{subfigure}
        \end{center}
    \caption{Wind usage amounts of the charging/discharging-based formulation under different pricing schemes} 
    \label{fig:wind_usage_stoch_charge_comb}
\end{figure}

We also analyze the effect of flat pricing in the charging/discharging-based 2SSP formulation and similarly observe less fluctuation in the generation and charging amounts. Specifically, we experimented with alternating prices in both of the two stages;  Figure~\ref{fig:wind_usage_stoch_charge_comb} gives an overview. 

It is immediate to see that wind usage profiles in Figure~\ref{fig:wind_usage_stoch_charge_comb} do not alter dramatically. Interestingly, when we calculate the wind utilization ratios explained in equation~\eqref{e:windUt}, all of the different pricing approaches result in the same utilization, 0.9645. This suggests that even though the usage of wind alters in shape, the total utilization remains the same. It is safe to conclude that the actual limitation of wind usage is rooted in the operational requirements of the transit fleet.

\begin{table}
\caption{Solution times of the charging/discharging-based formulation under different pricing schemes (s)}
\label{tab:cpu-time-charge-pricing}
\begin{center}
\begin{tabular}{l r  r r r}
    \toprule
    \hline
    Instance     & Flat both & Flat first & Flat second  & Peak both\\
    \hline 
     \texttt{case9}    &  101.4824 & 114.1081 & 114.0692 & 121.3390\\ 
     \texttt{case14}    & 27.0704 & 951.4918 & 22.0589 & 1690.1862\\ 
     \texttt{case30}    & 45.8860 & 40.0843 &41.3030  & 32.7924\\ 
     \texttt{case39}    & 63.4547  & 41.5623 & 42.6833 & 47.3546\\ 
     \texttt{case57}    & 1152.9587  & 166.6651 & 123.7377 & 52.4875\\ 
     \texttt{case118}     & 1443.3248  & 3216.8284 & 934.7387 & 991.3013\\ 
     \bottomrule
\end{tabular}
\end{center}
\end{table}

We present the solution times of the charging/discharging-based formulation under different pricing in Table~\ref{tab:cpu-time-charge-pricing}. One can observe that there is no single trend in solution times regardless of which pricing methodology is chosen. However, one should also note that some of the parameters chosen in these cases were problem-specific and could be the source of inconsistency. 

\subsection{Congestion analysis}
\label{sec:congestion-analysis}

\begin{figure}
\begin{center}
\includegraphics[trim = 0cm 0cm 0cm -0.2cm, width=10cm]{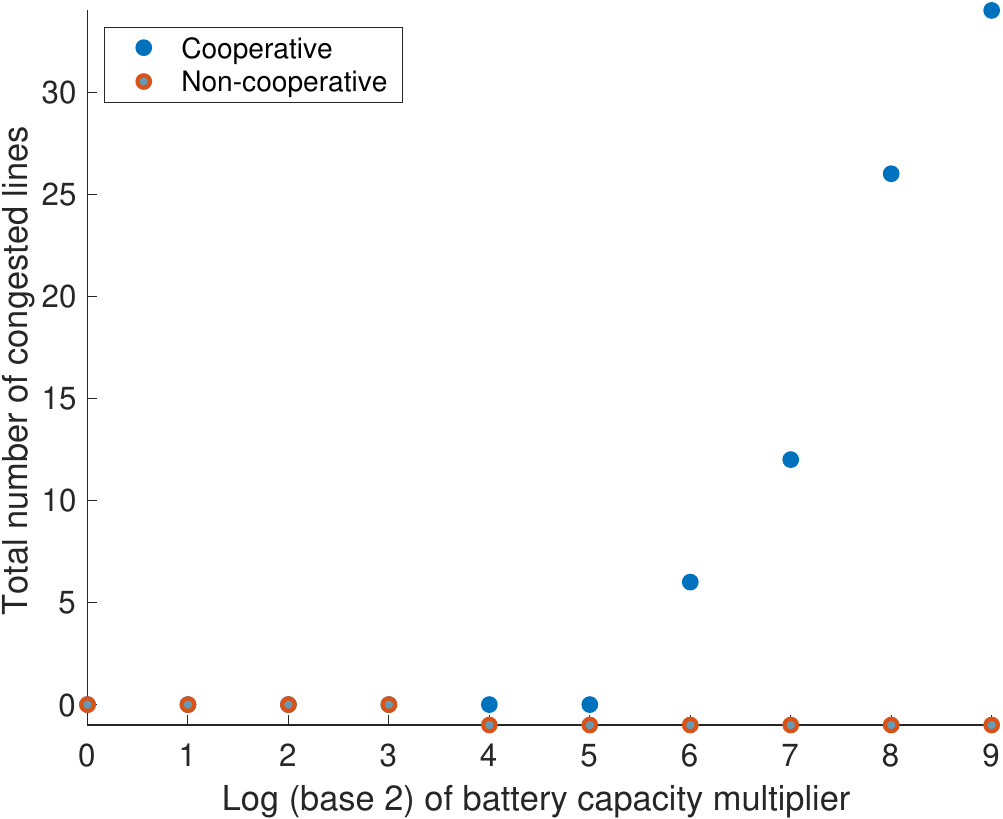}
\end{center}
\caption{Congestion analysis by varying battery capacity}
\label{fig:expansion_congestion}
\end{figure}

At present, both the small quantity of BEBs being utilized in the transit fleet and the limited battery capacity of these BEBs result in the BEBs not presenting any challenges to the grid by way of congestion. But the question is raised, as more BEBs are integrated into the transit fleet, and the battery technology continues to improve, at what point will the operation of BEBs begin to stress the power system? To this end, in our case files, we fix the line capacities at 10 MW and systematically increase the battery capacities of the current fleet, to represent a form of expansion, while holding power network parameters (and hence, capabilities) constant. Then, we solve the deterministic formulation and collect the total number of lines with full utilization over a day. Figure~\ref{fig:expansion_congestion} gives an illustration for 9-node instance.

We conduct this analysis for both cooperative and non-cooperative cases, by using the same prices in both. There are two main conclusions: firstly, the congestion in the network increases with larger battery capacities; and second, cooperation can mitigate the extra demand introduced since the non-cooperative model becomes infeasible after a certain point (where $-1$ indicates infeasibility). We have also investigated other power networks and obtained similar conclusions.

\section{Conclusion}\label{sec:con}
In this work, we propose a deterministic mathematical programming formulation to co-optimize the operation of the public transit system and the power grid when the electric transit buses are in their off-schedule considering generators without uncertainty. Furthermore, we propose two different two-stage stochastic programming formulations, in which uncertainty is present within the generation units with different recourse actions: \textit{(i)} ramping up/down the conventional generators, \textit{(ii)} additional charging/discharging of the transit fleet. 

Our numerical results demonstrate that coordinated operation between the ISO and the transit authority can decrease the generation costs while ensuring that the charging can be done at no detriment to the power system. Additionally, we conclude that in the presence of renewable generators, charging/discharging of the fleet as a recourse action can serve as a useful alternative to ramping up/down the conventional generators. Yet, as detailed, the benefits of charging/discharging of the fleet as a recourse action depend on battery capacities, ramping costs, and ramping limits. Our proposed formulations are general in the sense that they can incorporate different schedules, and different renewable generation units, different transportation, and power networks. 

We explore the effect of pricing the charging/discharging of the transit authority on the characteristics of the co-optimization for two different pricing strategies, namely peak pricing and flat pricing. Further, we also underline the effect of cooperation compared to a non-cooperative strategy by conducting a congestion analysis considering larger battery capacities. We found out that in the case of a non-cooperative strategy, the power dispatch problem becomes infeasible earlier in the electrification endeavor whereas in the cooperative approach, the fleet can charge/discharge at different points to ensure a feasible power flow with higher costs. 

A valuable extension to these formulations would be to consider separate renewable generation units owned and managed by the transportation authority. In this case, the charging/discharging costs of the vehicles would depend on which generator is chosen to be used. That is, it would be the same as the generation cost if they use their own generators, and it would be the market price derived from locational marginal prices if they use other generators within the power grid operated by the ISO.

Moreover, due to the fine resolution of the formulation, incorporating a physical model for battery degradation within the framework would decrease the losses further, and increase the flexibility of the co-optimization overall. Within the capabilities of the current framework, external effects such as battery degradation and charging technology could be represented by updating the necessary parameters. Finally, it is worth noting that while our approach does not address the decision of locating charging stations, the solutions obtained from our model would provide valuable managerial insights about the charging station locations.

\bibliography{coop}

\begin{appendices}

\section{Nomenclature}
\label{apx:nomenclature}

\begin{table}[H]
\caption{Nomenclature}
\begin{center}
\begin{tabular}{ l | l}
\hline
 \textbf{Sets} & \\ \hline 
 $\Ncal$ & set of nodes in the power network \\
 $\Ncal_g \subseteq \Ncal$ & set of nodes with conventional generators in the power network \\
 $\Ncal_r \subseteq \Ncal$ & set of nodes with renewable generators in the power network \\
 $\Ncal_b \subseteq \Ncal$ & set of candidate points for transit bus connection coupling the two systems \\ 
 $\Lcal$ & set of lines present in the power network \\
 $\Bcal$ & set of transit buses \\
 $\Tcal$ & set of time periods \\
 $\Tcal_b$ & set of off-schedule time periods for bus $b \in \Bcal$ \\
 $\Omega$ & set of scenarios for uncertain renewable generation \\ \hline
 \textbf{Indices} & \\ \hline
 $T_b^1, T_b^2$ & first and last time periods in off-schedule of bus $b \in \Bcal$, respectively \\ \hline
 \textbf{Parameters} & \\ \hline
 $c^g$ & conventional generation cost coefficients (linear)\\
 $c'^{g}$ & conventional generation cost coefficients (quadratic)\\
 $c$ & charging/discharging costs of transit buses\\
 $c^{+}$ & additional charging/discharging costs of transit buses\\
 $c^r$ & renewable generation cost coefficients \\
 $c^{\text{shed}}$ & costs of load-shedding in the power network\\
 $c^{g,+}, c^{g,-}$ & ramp-up and ramp-down costs of conventional generators, respectively\\
 $p^d$ & external demands present in the power grid\\
 $\overline{\theta}$ & limit on the voltage angles in the power network\\
 $x$ & reactance values of lines in the power network\\
 $\overline{S}$ & flow limits on lines in the power network\\
 $\overline{p}^g$ & conventional generation limits \\
 $p^{\delta g}$ & ramping limits of conventional generators \\
 $e^1$ & initial battery levels of transit buses (at the beginning of their off-schedules)\\
 $\underline{e}, \overline{e}$ & lower and upper limits on battery levels of transit buses, respectively\\
 $s$ & traversal energy consumption of transit buses over one time period \\
 $\overline{p}^c, \overline{p}^{dc}$ & limits on charging and discharging amounts of transit buses, respectively\\
 $\eta$ & efficiency values of transit bus batteries\\ 
 $\delta t$ & duration of time in a single period\\
 $\Delta t(\cdot,\cdot)$ & required number of time steps to relocate between two charging stations\\
 $\alpha$ & convex combination coefficient \\
 $\pi $ & probabilities of scenarios \\ \hline
 \textbf{Variables} & \\ \hline
 $p^g, p^r$ & conventional and renewable generation amounts, respectively \\
 $p^{d,\text{shed}}$ & amounts of load-shedding \\
 $p^{g,+}, p^{g,-}$ & ramp-up and ramp-down amounts of conventional generators, respectively\\
 $p^c, p^{dc}$ & charging and discharging amounts of transit buses, respectively\\
 $p^{c,+}, p^{dc,+}$ & additional charging and discharging amounts of transit buses, respectively\\
 $p$ & power flows on lines in the power network\\
 $\theta$ & voltage angles on nodes in the power network\\
 $e$ & battery levels of transit buses\\
 $z$ & assignment variables of transit buses to charging locations \\
 $y$ & indicator variables on the traversal status of transit buses \\ \hline
\end{tabular}
\end{center}
\end{table}

\section{Complete formulation of the deterministic model}

\label{apx:det_model}
This MIQP formulation is given by:
\begin{equation}
    \min \quad (1-\alpha)\sum_{t \in \Tcal} \sum_{i \in \Ncal_g} c^g_{it} p^g_{i t} + c'^g_{it} (p^g_{i t})^2 + \alpha \sum_{t \in \Tcal} \sum_{b \in \Bcal} \sum_{i \in \Ncal_b} c_{i t} \left(p^c_{ibt} - p^{dc}_{ibt} \right) \label{eq:obj1}
\end{equation}
subject to:
\begin{align}
    & p^g_{it} - \sum_{b \in \Bcal} p^c_{ibt} + \sum_{b \in \Bcal} p^{dc}_{ibt} - p^d_{it} = \sum_{j : (i,j) \in \mathcal{L}} p_{ijt} - \sum_{j : (j,i) \in \mathcal{L}} p_{jit}
    & \forall i \in \Ncal, t \in \Tcal \label{eq:c2}\\ 
     &p^g_{it} = 0  &\forall i \in \Ncal \setminus \Ncal_g, t \in \Tcal \label{eq:c3}\\
     & -\overline{\theta} \leq \theta_{it} \leq \overline{\theta} & \forall i \in \Ncal, t \in \Tcal \label{eq:c4}\\
     & \theta_{1t} = 0 &\forall t \in \Tcal \label{eq:c5}\\
     & p_{ijt} = \frac{\theta_{it} - \theta_{jt}}{x_{ij}} & \forall (i,j) \in \Lcal, t \in \Tcal \label{eq:c6}\\
    & -\overline{S}_{ij} \leq p_{ijt} \leq \overline{S}_{ij} & \forall (i,j) \in \Lcal, t \in \Tcal \label{eq:c7}\\
    &0 \leq p^g_{it} \leq \overline{p}_{it}^g &\forall i \in \Ncal_g, t \in \Tcal \label{eq:c8}\\
    & -p^{\delta g}_{it} \leq p^g_{i,t+1} - p^g_{it}\leq p^{\delta g}_{it} &\forall i \in \Ncal_g, t \in \Tcal \setminus \{T\} \label{eq:c9}\\
     &e_{b T^1_b} = e^1_b &\forall b \in \Bcal \label{eq:c10}\\
     &e_{b T^2_b} + \eta_b \sum_{i \in \Ncal_b} p^c_{ib T^2_b} \delta t - \frac{1}{\eta_b} \sum_{i \in \Ncal_b} p^{dc}_{ibT^2_b} \delta t - s_b  y_{bT^2_b} = \overline{e}_b &\forall b \in \Bcal \label{eq:c11}\\
    &e_{b,t+1} = e_{bt} + \eta_b \sum_{i \in \Ncal_b} p^c_{ibt} \delta t - \frac{1}{\eta_b} \sum_{i \in \Ncal_b} p^{dc}_{ibt} \delta t - s_b  y_{bt} &\forall b \in \Bcal, \: t,t+1 \in \Tcal_b \label{eq:c12}\\
    &\underline{e}_{b} \leq e_{bt} \leq \overline{e}_{b}  &\forall b \in \Bcal, t \in \Tcal_b \label{eq:c13}\\
    &0 \leq  p^c_{ibt} \leq \overline{p}^c_{b} z_{ibt} &\forall i \in \Ncal_b, b \in \Bcal, t \in \Tcal \label{eq:c14}\\
    &0 \leq  p^{dc}_{ibt} \leq \overline{p}^{dc}_{b} z_{ibt} &\forall i \in \Ncal_b, b \in \Bcal, t \in \Tcal \label{eq:c15}\\
     &p^c_{ibt} = 0, ~ p^{dc}_{ibt} = 0 &\forall i \in \Ncal \setminus \Ncal_b, b \in \Bcal, t \in \Tcal \label{eq:c16}\\
     &\sum_{i \in \Ncal_b} z_{ibt} + y_{bt} = 1&\forall b \in \Bcal, t \in \Tcal_b \label{eq:c17}\\
     &z_{ibt} + z_{jbt'} \leq 1 &  \forall t' \in \Tcal_b, ~ t < t' \leq t+\Delta t(i,j), \nonumber \\
     & & \forall i,j \in \Ncal_b, i \neq j,  b \in \Bcal, t \in \Tcal_b \label{eq:c18} \\
     & z_{d b T^1_b} = 1 &  d = \Ncal_{b}(1), \forall b \in \Bcal \label{eq:c19}\\
     &y_{bt} = 0 &\forall b \in \Bcal, t \in \Tcal \setminus \Tcal_b \label{eq:c20}\\
     &z_{ibt} = 0 &\forall i \in \Ncal_b, b \in \Bcal, t \in \Tcal \setminus \Tcal_b \label{eq:c21}\\
     &y_{bt} \in \{0,1\} & \forall b \in \Bcal, t \in \Tcal \label{eq:c22}\\
     &z_{ibt} \in \{0,1\} & \forall i \in \Ncal_b, b \in \Bcal, t \in \Tcal \label{eq:c23}
\end{align}

Our objective \eqref{eq:obj1} is to minimize a convex combination, with coefficient $\alpha$, of the total power generation cost and the charging/discharging cost of the transit buses. For each time period $t \in \Tcal$, both the charging and discharging costs are given by $c_{it}$. These values, also known as locational marginal prices \citep{lamadrid-ancillary}, represent the optimal dual values associated with the nodal balance equations of the power network. Hence, in the case of simultaneous charging and discharging, the costs and revenues will cancel each other out; that is, it is never (strictly) optimal to charge and discharge simultaneously. This simplifies the formulation since otherwise, one would need an extra set of binary variables indicating whether each vehicle is charging or discharging.

Constraints \eqref{eq:c2}-\eqref{eq:c8} are standard DC optimal power flow constraints except constraint set \eqref{eq:c2} additionally incorporates terms for charging and discharging in the nodal balance equations. Inequality \eqref{eq:c9} ensures that the ramping amount of generators is within the limit. 
Constraint sets \eqref{eq:c10} and \eqref{eq:c11} provide initial and final conditions on the battery levels of the transit buses respectively. Constraints \eqref{eq:c12} are battery level updates for the transit fleet, where the bounds on the battery levels are employed in \eqref{eq:c13}. Constraints \eqref{eq:c14} and \eqref{eq:c15} ensure that a battery can only charge/discharge if the transit bus is connected to a node in the power network. 

Constraints \eqref{eq:c17} are assignment constraints specifying that a transit bus can either be connected to one of the possible nodes in the power network or relocating within the network. Then, constraint set \eqref{eq:c18} guarantees any relocation made throughout the horizon is feasible while ensuring that a transit bus cannot relocate in fewer time steps than the required travel time $\Delta t(\cdot,\cdot)$. Finally, equation \eqref{eq:c19} is the initial assignment of the transit buses to the depot node, whereas \eqref{eq:c20} and \eqref{eq:c21} ensure that assignment only occurs within the off-schedule periods of the vehicles. Note that the periods in the formulation are inherently cyclic (i.e. considering hourly intervals, period 23 connects to period 0).

As a note related to the assumption on the capacities of charging stations, considering constant charging capacities per station $C_i$, depending on the specific application, one can incorporate the following constraint set:
\begin{equation}
    \sum_{b \in \Bcal} z_{ibt} \leq C_i \hspace{2cm} \forall i \in \Ncal_b, t \in \Tcal
\end{equation}

\section{Complete formulation of the two-stage stochastic ramping-based formulation}
\label{apx:stoch_ramping_model}

The complete formulation is as follows:
\begin{align}
    \min \quad & (1-\alpha) \left [ \sum_{t \in \Tcal}  \sum_{i \in \Ncal_g} c^g_{it} p^g_{i t} + c'^g_{it} (p^g_{i t})^2 + \sum_{\omega \in \Omega} \pi_\omega \left ( \sum_{t \in \Tcal} \sum_{i \in \Ncal_r} c^r_{it} p^r_{it\omega}  +  \sum_{t \in \Tcal}\sum_{i \in \Ncal }c_{it}^{\text{shed}} p_{it\omega}^{d,\text{shed}} \right ) \right ] \nonumber \\
    & +  (1-\alpha) \sum_{\omega \in \Omega }\pi_\omega  \sum_{t \in \Tcal} \sum_{i \in \Ncal_g} \left( c_{i t}^{g,+} p^{g,+}_{it\omega} - c_{i t}^{g,-} p^{g,-}_{it\omega}    \right ) + \alpha \sum_{t \in \Tcal} \sum_{b \in \Bcal} \sum_{i \in \Ncal_b} c_{i t} \left(p^c_{ibt} - p^{dc}_{ibt} \right) \label{eq:obj1_s1}
\end{align}
subject to:
\begin{align}
& p^g_{it} + p^r_{it} - \sum_{b \in \Bcal} p^c_{ibt} + \sum_{b \in \Bcal} p^{dc}_{ibt} - p^d_{it} = \sum_{j : (i,j) \in \mathcal{L}} p_{ijt} - \sum_{j : (j,i) \in \mathcal{L}} p_{jit}
    & \forall i \in \Ncal, t \in \Tcal \label{eq:sc1}\\  
    & p^r_{it\omega}-p^r_{it} + p^{g,+}_{it\omega} -  p^{g,-}_{it\omega} + p^{d,shed}_{it\omega}= \sum_{j : (i,j) \in \mathcal{L}} \left (p_{ijt\omega} - p_{ijt} \right ) - \sum_{j : (j,i) \in \mathcal{L}} \left ( p_{jit\omega} - p_{jit}\right )
    & \forall i \in \Ncal, t \in \Tcal, \omega \in \Omega \label{eq:sc2}
    \end{align}
    \begin{align}
     &p^g_{it} = 0  &\forall i \in \Ncal \setminus \Ncal_g, t \in \Tcal \label{eq:sc3}\\
     &p^r_{it} = 0  &\forall i \in \Ncal \setminus \Ncal_r, t \in \Tcal \label{eq:sc4}\\
     &p^r_{it\omega} = 0  &\forall i \in \Ncal \setminus \Ncal_r, t \in \Tcal, \omega \in \Omega \label{eq:sc5}\\
     & -\overline{\theta} \leq \theta_{it} \leq \overline{\theta} & \forall i \in \Ncal, t \in \Tcal \label{eq:sc6}\\
     & -\overline{\theta} \leq \theta_{itw} \leq \overline{\theta} & \forall i \in \Ncal, t \in \Tcal, \omega \in \Omega \label{eq:sc7}\\
     & \theta_{1t} = 0 &\forall t \in \Tcal \label{eq:sc8}\\
     & \theta_{1t\omega} = 0 &\forall t \in \Tcal, \omega \in \Omega \label{eq:sc9}\\
     & p_{ijt} = \frac{\theta_{it} - \theta_{jt}}{x_{ij}} & \forall (i,j) \in \Lcal, t \in \Tcal \label{eq:sc10}\\
       & p_{ijt\omega} = \frac{\theta_{it\omega} - \theta_{jt\omega}}{x_{ij}} & \forall (i,j) \in \Lcal, t \in \Tcal, \omega \in \Omega \label{eq:sc11}\\
    & -\overline{S}_{ij} \leq p_{ijt} \leq \overline{S}_{ij} & \forall (i,j) \in \Lcal, t \in \Tcal \label{eq:sc12}\\
    & -\overline{S}_{ij} \leq p_{ijt\omega} \leq \overline{S}_{ij} & \forall (i,j) \in \Lcal, t \in \Tcal, \omega \in \Omega \label{eq:sc13}\\
    &0 \leq p^g_{it} + p^{g,+}_{it} \leq \overline{p}_{it}^g &\forall i \in \Ncal_g, t \in \Tcal \label{eq:sc14}\\
    & p^{g,-}_{it} \leq p^g_{it} &\forall i \in \Ncal_g, t \in \Tcal \label{eq:sc15}\\
      &  p^{g,+}_{it\omega} \leq p^{g,+}_{it} &\forall i \in \Ncal_g, t \in \Tcal, \omega \in \Omega \label{eq:sc16}\\
    &  p^{g,-}_{it\omega} \leq p^{g,-}_{it} &\forall i \in \Ncal_g, t \in \Tcal, \omega \in \Omega \label{eq:sc17}\\
    & -p^{\delta g}_{it} \leq p^g_{i,t+1} - p^g_{it}\leq p^{\delta g}_{it} &\forall i \in \Ncal_g, t \in \Tcal \setminus \{T\} \label{eq:sc18}\\
    &0 \leq p^r_{it} \leq \overline{p}_{it}^r &\forall i \in \Ncal_r, t \in \Tcal \label{eq:sc19}\\
    &0 \leq p^r_{it\omega} \leq \tilde{p}_{it\omega}^r &\forall i \in \Ncal_r, t \in \Tcal, \omega \in \Omega \label{eq:sc20}\\
    &0 \leq p^{d,\text{shed}}_{it \omega} \leq p^d_{it} & \forall i \in \Ncal, t \in \Tcal, \omega \in \Omega \label{eq:sc21}\\
    &0 \leq p^{g,+}_{it} \leq \overline{p}^{g,+}_{it}  &\forall i \in \Ncal_g,  t \in \Tcal \label{eq:sc22}\\
    &0 \leq  p^{g,-}_{it} \leq  \overline{p}^{g,-}_{it} &\forall i \in \Ncal_g,  t \in \Tcal \label{eq:sc23}\\
     &p^{g,+}_{it\omega} \geq 0, ~ p^{g,-}_{it\omega} \geq 0 &\forall i \in \Ncal_g, t \in \Tcal, \omega \in \Omega \label{eq:sc24} \\
     &e_{b T^1_b} = e^1_b &\forall b \in \Bcal \label{eq:sc25}\\
     &e_{b T^2_b} + \eta_b \sum_{i \in \Ncal_b} p^c_{ib T^2_b} \delta t - \frac{1}{\eta_b} \sum_{i \in \Ncal_b} p^{dc}_{ibT^2_b} \delta t - s_b  y_{bT^2_b} = \overline{e}_b &\forall b \in \Bcal \label{eq:sc26}\\
    &e_{b,t+1} = e_{bt} + \eta_b \sum_{i \in \Ncal_b} p^c_{ibt} \delta t - \frac{1}{\eta_b} \sum_{i \in \Ncal_b} p^{dc}_{ibt} \delta t - s_b  y_{bt} &\forall b \in \Bcal, \: t,t+1 \in \Tcal_b \label{eq:sc27}\\
    &\underline{e}_{b} \leq e_{bt} \leq \overline{e}_{b}  &\forall b \in \Bcal, t \in \Tcal_b \label{eq:sc28}\\
    &0 \leq  p^c_{ibt} \leq \overline{p}^c_{b} z_{ibt} &\forall i \in \Ncal_b, b \in \Bcal, t \in \Tcal \label{eq:sc29}\\
    &0 \leq  p^{dc}_{ibt} \leq \overline{p}^{dc}_{b} z_{ibt} &\forall i \in \Ncal_b, b \in \Bcal, t \in \Tcal \label{eq:sc30}\\
     &p^c_{ibt} = 0, ~ p^{dc}_{ibt} = 0 &\forall i \in \Ncal \setminus \Ncal_b, b \in \Bcal, t \in \Tcal \label{eq:sc31}\\
     &\sum_{i \in \Ncal_b} z_{ibt} + y_{bt} = 1&\forall b \in \Bcal, t \in \Tcal_b \label{eq:sc32}\\
     &z_{ibt} + z_{jbt'} \leq 1 &  \forall t' \in \Tcal_b, ~ t < t' \leq t+\Delta t(i,j), \label{eq:sc33}\\
     & & \forall i,j \in \Ncal_b, i \neq j,  b \in \Bcal, t \in \Tcal_b \nonumber \\
     & z_{d b T^1_b} = 1 &  d = \Ncal_{b}(1), \forall b \in \Bcal \label{eq:sc34}\\
     &y_{bt} = 0 &\forall b \in \Bcal, t \in \Tcal \setminus \Tcal_b \label{eq:sc35}\\
     &z_{ibt} = 0 &\forall i \in \Ncal_b, b \in \Bcal, t \in \Tcal \setminus \Tcal_b \label{eq:sc36}\\
     &y_{bt} \in \{0,1\} & \forall b \in \Bcal, t \in \Tcal \label{eq:sc37}\\
     &z_{ibt} \in \{0,1\} & \forall i \in \Ncal_b, b \in \Bcal, t \in \Tcal \label{eq:sc38}
\end{align}

Here, the objective function \eqref{eq:obj1_s1} is the summation of first-stage costs (which has already been captured by the deterministic objective function \eqref{eq:obj1}) and second-stage costs including the expected renewable generation costs and expected ramping costs. In addition to the constraints captured in the deterministic model in Section~\ref{sec:model_det} are the second-stage nodal balance \eqref{eq:sc2} and flows \eqref{eq:sc7}, \eqref{eq:sc9}, \eqref{eq:sc11}, \eqref{eq:sc13}; limits on renewable generation \eqref{eq:sc19}, \eqref{eq:sc20}; limits on second-stage ramping \eqref{eq:sc14}, \eqref{eq:sc16}, \eqref{eq:sc17}; and limits on shedding in the second stage \eqref{eq:sc21}. The main differences from the two-stage stochastic OPF formulation presented in \cite{morales2013integrating} are constraints on the transit bus battery levels \eqref{eq:sc25}-\eqref{eq:sc28}, charging limits \eqref{eq:sc29}-\eqref{eq:sc31}, and transportation constraints \eqref{eq:sc32}-\eqref{eq:sc38}. More importantly, we have each of the OPF constraints repeated for multiple time periods rather than a single time period presented by \cite{morales2013integrating}. The inclusion of multiple periods, and their coupling by the presence of batteries, increase the complexity of the problem dramatically.

Note that we have two sets of ramping quantities $p_{it}^{\delta g}$, and variables $(p_{it}^{g,+}, p_{it}^{g,-})$ associated with ramping in different time periods. Specifically, the parameters $p_{it}^{\delta g}$ serve as upper bounds on the ramping of conventional generators between two time periods in the formulation, whereas the variables $(p_{it}^{g,+}, p_{it}^{g,-})$ correspond to ramping of conventional generators between the two stages of the stochastic formulation.

\section{Two-stage stochastic multi-period OPF formulation to obtain prices of charging/discharging}

\label{apx:MPOPF}

In this section, we present a baseline OPF formulation in order to estimate the first-stage and second-stage prices of charging/discharging in formulations presented in Section~4. The complete formulation is as follows:
\begin{align}
    \min \quad & \sum_{t \in \Tcal} \sum_{i \in \Ncal_g} c^g_{it} p^g_{i t} + c'^g_{it} (p^g_{i t})^2 + \sum_{\omega \in \Omega} \pi_\omega \left ( \sum_{t \in \Tcal} \sum_{i \in \Ncal_r} c^r_{it} p^r_{it\omega}  +  \sum_{t \in \Tcal}\sum_{i \in \Ncal }c_{it}^{\text{shed}} p_{it\omega}^{d,\text{shed}} \right )  \nonumber \\
    & +  \sum_{\omega \in \Omega }\pi_\omega  \sum_{t \in \Tcal} \sum_{i \in \Ncal_g} \left( c_{i t}^{g,+} p^{g,+}_{it\omega} - c_{i t}^{g,-} p^{g,-}_{it\omega}    \right ) \label{eq:obj1_s}
\end{align}
subject to:
\begin{equation}
     (c_{it}): ~~  p^g_{it} + p^r_{it} - p^d_{it} = \sum_{j : (i,j) \in \mathcal{L}} p_{ijt} - \sum_{j : (j,i) \in \mathcal{L}} p_{jit} 
     \hspace{2.9cm} \forall i \in \Ncal, t \in \Tcal \label{eq:nodal_balance_1}
\end{equation}
\begin{multline}
     (c_{it\omega}^+): ~~ p^r_{it\omega}-p^r_{it} + p^{g,+}_{it\omega} -  p^{g,-}_{it\omega} + p^{d,shed}_{it\omega}= \sum_{j : (i,j) \in \mathcal{L}} \left (p_{ijt\omega} - p_{ijt} \right ) - \sum_{j : (j,i) \in \mathcal{L}} \left ( p_{jit\omega} - p_{jit}\right )
    \\ \forall i \in \Ncal, t \in \Tcal, \omega \in \Omega  \label{eq:nodal_balance_2}
\end{multline}
    \begin{align}
     &p^g_{it} = 0  &\forall i \in \Ncal \setminus \Ncal_g, t \in \Tcal \\
     &p^r_{it} = 0  &\forall i \in \Ncal \setminus \Ncal_r, t \in \Tcal \\
     &p^r_{it\omega} = 0  &\forall i \in \Ncal \setminus \Ncal_r, t \in \Tcal, \omega \in \Omega \\
     & -\overline{\theta} \leq \theta_{it} \leq \overline{\theta} & \forall i \in \Ncal, t \in \Tcal \\
     & -\overline{\theta} \leq \theta_{itw} \leq \overline{\theta} & \forall i \in \Ncal, t \in \Tcal, \omega \in \Omega \\
     & \theta_{1t} = 0 &\forall t \in \Tcal \\
     & \theta_{1t\omega} = 0 &\forall t \in \Tcal, \omega \in \Omega\\
     & p_{ijt} = \frac{\theta_{it} - \theta_{jt}}{x_{ij}} & \forall (i,j) \in \Lcal, t \in \Tcal \\
       & p_{ijt\omega} = \frac{\theta_{it\omega} - \theta_{jt\omega}}{x_{ij}} & \forall (i,j) \in \Lcal, t \in \Tcal, \omega \in \Omega \\
    & -\overline{S}_{ij} \leq p_{ijt} \leq \overline{S}_{ij} & \forall (i,j) \in \Lcal, t \in \Tcal \\
    & -\overline{S}_{ij} \leq p_{ijt\omega} \leq \overline{S}_{ij} & \forall (i,j) \in \Lcal, t \in \Tcal, \omega \in \Omega \\
    &0 \leq p^g_{it} + p^{g,+}_{it} \leq \overline{p}_{it}^g &\forall i \in \Ncal_g, t \in \Tcal \\
    & p^{g,-}_{it} \leq p^g_{it} &\forall i \in \Ncal_g, t \in \Tcal \\
      &  p^{g,+}_{it\omega} \leq p^{g,+}_{it} &\forall i \in \Ncal_g, t \in \Tcal, \omega \in \Omega \\
    &  p^{g,-}_{it\omega} \leq p^{g,-}_{it} &\forall i \in \Ncal_g, t \in \Tcal, \omega \in \Omega \\
    & -p^{\delta g}_{it} \leq p^g_{i,t+1} - p^g_{it}\leq p^{\delta g}_{it} &\forall i \in \Ncal_g, t \in \Tcal \setminus \{T\} \\
    &0 \leq p^r_{it} \leq \overline{p}_{it}^r &\forall i \in \Ncal_r, t \in \Tcal \\
    &0 \leq p^r_{it\omega} \leq \tilde{p}_{it\omega}^r &\forall i \in \Ncal_r, t \in \Tcal, \omega \in \Omega \\
    &0 \leq p^{d,\text{shed}}_{it \omega} \leq p^d_{it} & \forall i \in \Ncal, t \in \Tcal, \omega \in \Omega\\
    &0 \leq p^{g,+}_{it} \leq \overline{p}^{g,+}_{it}  &\forall i \in \Ncal_g,  t \in \Tcal \\
    &0 \leq  p^{g,-}_{it} \leq  \overline{p}^{g,-}_{it} &\forall i \in \Ncal_g,  t \in \Tcal \\
     &p^{g,+}_{it\omega} \geq 0, ~ p^{g,-}_{it\omega} \geq 0 &\forall i \in \Ncal_g, t \in \Tcal, \omega \in \Omega
\end{align}

This formulation is a direct extension of the two-stage single-period optimal power formulation presented by \cite{morales2013integrating}. Specifically, we use the optimal dual variables $c_{it}, c_{it\omega}^+$ associated with constraints \eqref{eq:nodal_balance_1} and \eqref{eq:nodal_balance_2}, respectively, to determine the first-stage and second-stage charging/discharging prices.

\end{appendices}

\end{document}